\documentclass[11pt]{article}
\pdfoutput=1
\usepackage{amsmath,amssymb}

\usepackage{graphicx}
\usepackage{dcolumn}
\usepackage{bm}
\usepackage{epstopdf,overpic,color,rotating}

\graphicspath{{figures/}}

\newcommand*\samethanks[1][\value{footnote}]{\footnotemark[#1]}

\title{Compressive sensing and low-rank libraries for classification \\ of bifurcation regimes in nonlinear dynamical systems}

\author{Steven L. Brunton$^\ddagger$\thanks{Department of Applied Mathematics, University of Washington, Seattle, WA. 98195-2420.  $^\ddagger$ ({sbrunton@uw.edu}). Questions, comments, or corrections to this document may be directed to that email address.} \and Jonathan H. Tu\thanks{Mechanical and Aerospace Engineering, Princeton University, Princeton, NJ 08544 } \and Ido Bright\samethanks[1] \and J. Nathan Kutz\samethanks[1]}

\date{December 14, 2013}

\begin{document}
\maketitle

\begin{abstract}
We show that for complex nonlinear systems, model reduction and compressive sensing strategies can be combined to great advantage for classifying, projecting, and reconstructing the relevant low-dimensional dynamics.  $\ell_2$-based dimensionality reduction methods such as the proper orthogonal decomposition are used to construct separate modal libraries and Galerkin models based on data from a number of bifurcation regimes.  These libraries are then concatenated into an over-complete library, and $\ell_1$ sparse representation in this library from a few noisy measurements results in correct identification of the bifurcation regime.  This technique provides an objective and general framework for classifying the bifurcation parameters, and therefore, the underlying dynamics and stability.  After classifying the bifurcation regime, it is possible to employ a low-dimensional Galerkin model, only on modes relevant to that bifurcation value.  These methods are demonstrated on the complex Ginzburg-Landau equation using sparse, noisy measurements.  In particular, three noisy measurements are used to accurately classify and reconstruct the dynamics associated with six distinct bifurcation regimes; in contrast, classification based on least-squares fitting ($\ell_2$) fails consistently.
\end{abstract}


\section{Introduction}
Nonlinear dynamical systems are ubiquitous in characterizing the behavior of physical, biological and engineering systems.
With few exceptions, nonlinearity impairs our ability to construct analytically tractable solutions, and we instead rely on experiments and high-performance computation to study a given system.
Numerical discretization can often yield a system of equations with millions or billions of degrees of freedom.  Thus, both simulations and experiments can generate enormous data sets that strain computational resources and confound one's understanding of the underlying dynamics.
Fortunately, many high-dimensional systems exhibit dynamics that evolve on a slow-manifold and/or a low-dimensional attractor (e.g., pattern forming systems~\cite{Cross:1993}).  
We propose a data-driven modeling strategy that represents low-dimensional dynamics using {\em dimensionality reduction} methods such as the proper orthogonal decomposition (POD)~\cite{HLBR_turb} and classifies/reconstructs the observed low-dimensional manifolds with {\em compressive (sparse) sensing} (CS)~\cite{Donoho:2006,Candes:2006c,Baraniuk:2009,Tropp:2007}.
Thus dynamic structures are represented efficiently with the $\ell_2$ norm and identified from sparse measurements with the $\ell_1$ norm.

The application of machine-learning and compressive sensing to dynamical systems is synergistic, in that underlying low-rank structures facilitate sparse measurements~\cite{Kutz:2013}.  This combination has the potential to transform a number of challenging fields.  
Such a strategy may enhance nonlinear estimation and control, where real-time analysis is critical.  
Moreover, adaptive time-stepping algorithms can take advantage of the low-dimensional embedding for greatly reduced computational costs~\cite{Kevrekidis:2003,Schaeffer:2013}.  
Additionally, the interplay of sparsity and complex systems has been investigated with the goal of overcoming the curse of dimensionality associated with neuronal activity and neuro-sensory systems~\cite{Ganguli:2012}.  
Compressive sensing may also play a role in similar statistical learning, library-based, and/or information theory methods~\cite{Coifman:1992,Branicki:2013} used in fluid dynamics~\cite{Bright:2013,Glauser:2013}, climate science~\cite{Giannakis:2013,Branicki:2013} and oceanography~\cite{Agarwal:2011}.  
Indeed, compressive sensing is already playing a critical role in model building and assessment in the physical sciences~\cite{Nelson:2013,Shabani:2011,Wang:2011,Tayler:2012}.
These challenging open problems would benefit from a paradigm shift in modeling and analysis, whereby low-dimensional coherence is leveraged for use with sparse sampling techniques.

\subsection{Challenges of POD-Galerkin models across parameter regimes}  
Galerkin-POD is a well-known~\cite{HLBR_turb} dimensionality reduction method for complex systems.  In the context of fluid dynamics, Galerkin projection of the Navier-Stokes equations onto a truncated POD mode basis is an effective method of model-order reduction, resulting in a system of ordinary differential equations.  However, Galerkin projection onto POD modes obtained across a range of parameter values, the so-called global POD~\cite{Taylor:2002,Schmidt:2004}, often results in unstable and/or inaccurate models.  There have been a number of modifications to POD-Galerkin models that seek to address this issue, but it remains a major challenge of low-order modeling in fluids.  

A modified method that uses interpolated angles of multiple POD subspaces has been demonstrated to capture F-16 parameterized dynamics~\cite{Lieu:2005}.  Including additional modes, such as the shift mode~\cite{noack:03cyl}, to capture transients between qualitatively different flow regimes has resulted in additional methods such as double POD~\cite{Siegel:2008}, and the GaussÐNewton with approximated tensors (GNAT) method~\cite{Carlberg:2013} and trust-region POD~\cite{Fahl:2000,Bergmann:2008}.  Alternative methods for stabilizing POD by adding additional modes and closure terms have been investigated~\cite{Balajewicz:2012b,Osth:2014}.  In each case, the objective is to construct a dimensionally reduced set of dynamics that accurately represent the underlying complex system and that does not suffer from instabilities.

\subsection{Current approach}
To avoid a single POD-Galerkin model defined across dynamical regions, we instead develop a classification scheme to determine which dynamic region our system is in, and then use a Galerkin model defined only on modes in that region.  
The procedure advocated here involves two main steps.  
First, a modal library is constructed that is representative of a number of distinct dynamical regimes, i.e. the low-dimensional attractors are approximated by their optimal bases.  
Second, compressive sensing techniques are applied using this learned library.  
The goals are threefold: 1) classify the dynamic regime, 2) project the measurements onto the correct modal amplitudes, and 3) reconstruct the low-dimensional dynamics through Galerkin projection~\cite{HLBR_turb}.  
Here we concatenate POD bases to construct the library, although generalizing the library building strategy into a broader machine learning context~\cite{duda:pattern} is interesting and may yield even more efficient strategies.  
There are many ways to build a library, especially considering the three goals above.  
In this case, we keep distinct POD bases for each dynamic regime, since this is better for the Galerkin projection step.  The classification scheme, using $\ell_1$ minimization in a over-complete library, is closely related to sparse representation from image classification~\cite{Wright:2009}.

The paper is outlined as follows:  In Sec.~\ref{sec:background} a brief review is provided of the compressive sensing architecture and its relationship to $\ell_1$ convex optimization.  
Also reviewed are the basic ideas behind the proper orthogonal decomposition (POD) for $\ell_2$ dimensionality reduction.  
These methods are combined in Sec.~\ref{sec:methods} to form the key contributions of this work.  
Namely, the $\ell_2$-norm provides the sparse basis modes used by the $\ell_1$-norm for sparse representation.  
Section~\ref{sec:results} demonstrates the use of these techniques on one of the classical models of mathematical physics:  the Ginzburg-Landau equation.
An outlook of the advantages and general applicability of the method to complex systems is given in the concluding section~\ref{sec:discussion}.

\section{Background}\label{sec:background}
In the following subsections, we introduce two well-established techniques that will be combined in this paper.  The first method is compressive sensing (CS), whereby a signal that is sparse in some basis may be recovered using proportionally few measurements by solving for the $\ell_1$-minimizing solution to an underdetermined system.  The second method is the proper orthogonal decomposition (POD), which allows a dataset to be reduced optimally in an $\ell_2$ sense.  

Both theories have been applied to a range of problems.  In this paper, we advocate combining these methods since the $\ell_2$ basis obtained from POD is a particularly good choice of a sparse basis for compressive sensing.  The underlying reason for this is that the data is obtained from the low-dimensional attractors of the governing
complex system.

\subsection{$\ell_1$-based sparse sensing}\label{sec:back:l1}
Consider a high-dimensional measurement vector ${\bf U}\in\mathbb{R}^n$, which is sparse in some space, spanned by the columns of a matrix ${\bf \Psi}$:
\begin{equation}
{\bf U}(x,t) = {\bf \Psi a}.
\label{eq:sparse1}
\end{equation}
Here, sparsity means that $\bf U$ may be represented in the transform basis $\bf\Psi$ by a vector of coefficients $\bf a$ that contains mostly zeros.  
More specifically, $K$-\textit{sparsity} means that there are $K$ nonzero elements.  
In this sense, sparsity implies that the signal is \emph{compressible}.      

Consider a sparse measurement ${\bf \hat U}\in\mathbb{R}^m$, with $m\ll n$:
\begin{equation}
{\bf \hat U} = {\bf \Phi U},
\label{eq:sparse2}
\end{equation}
where $\Phi$ is a measurement matrix that maps the full state
measurement ${\bf U}$ to the sparse measurement vector ${\bf \hat U}$.  Details
of this measurement matrix will be given shortly.
Plugging \eqref{eq:sparse1} into \eqref{eq:sparse2} yields an underdetermined system:
\begin{equation}
{\bf \hat U} = {\bf \Phi\Psi a}.
\label{eq:sparse3}
\end{equation}

We may then solve for the sparsest solution ${\bf a}$ to the underdetermined system of equations in \eqref{eq:sparse3}.  Sparsity is measured by the $\ell_0$ norm, and solving for the solution ${\bf a}$ that has the smallest $|{\bf a}|_0$ norm is a combinatorially hard problem.  However, this problem may be relaxed to a convex problem, whereby the $|{\bf a}|_1$ norm is minimized, which may be solved in polynomial time~\cite{Candes:2006c,Donoho:2006}.  The specific minimization problem is:
\begin{equation*}
\arg \min {|{\bf \hat a}|_1} \text{ such that } {\bf \Phi\Psi\hat a}={\bf \hat U}.
\end{equation*}
There are other algorithms that result in sparse solution vectors, such as orthogonal matching pursuit~\cite{Tropp:2007}.

This procedure, known as {\em compressive sensing}, is a recent development that has had widespread success across a range of problems.  There are technical issues that must be addressed.  For example, the number of measurements $m$ in ${\bf \hat U}$ should be on the order of $K\log(n/K)$, where $K$ is the degree of sparsity of $\bf a$ in $\bf \Psi$~\cite{Candes:2006, Candes:2006a,Baraniuk:2007}.  In addition, the measurement matrix $\bf \Phi$ must be {\em incoherent} with respect to the sparse basis $\bf \Psi$, meaning that the columns of $\bf \Phi$ and the columns of $\bf \Psi$ are uncorrelated.  
Interestingly, significant work has gone into demonstrating that Bernouli and Gaussian random measurement matrices are almost certainly incoherent with respect to a given basis~\cite{Candes:2006b}. 

Typically a generic basis such as Fourier or wavelets is used in conjunction with sparse measurements consisting of random projections of the state.  However, in many engineering applications, it is unclear how random projections may be obtained without first starting with a dense measurement of the state.  In this work, we constrain the measurements to be point measurements of the state, so that $\bf \Phi$ consists of rows of a permutation matrix.  
Our primary motivation for such point measurements arises from physical considerations in such applications as ocean or atmospheric monitoring where point measurements are physically relevant.  Moreover, sparse sensing is highly desirable as each measurement device is often prohibitively expensive, thus motivating much of our efforts in using sparse measurements to characterize the complex dynamics.

\subsection{$\ell_2$-based dimensionality reduction}
The proper orthogonal decomposition (POD)~\cite{Lumley:1970,HLBR_turb} is a tool with ubiquitous use in dimensionality reduction of physical systems\footnote{POD is sometimes referred to as principal components analysis~\cite{Pearson:1901}, the Karhunen--Lo\`eve decomposition, empirical orthogonal functions~\cite{lorenzMITTR56}, or the Hotelling transform~\cite{hotellingJEdPsy33_1,hotellingJEdPsy33_2}.}.  Data snapshots ${\bf U}(x,t_1), {\bf U}(x,t_2),\cdots {\bf U}(x,t_q)$ are collected into columns of a matrix ${\bf A}\in\mathbb{R}^{n\times q}$.  We then compute the singular value decomposition (SVD) of $\bf A$:
\begin{equation*}
{\bf A} = {\bf\Psi \Sigma W}^*.
\end{equation*}
Columns of the matrix ${\bf \Psi}$ are POD modes\footnote{Often, POD modes are given by the matrix ${\bf \Phi}$.  However, we choose $\bf \Psi$ for the POD basis and $\bf \Phi$ for the sparse measurement matrix for consistency with compressive sensing literature.  This is not to be confused with notation from balanced POD, where ${\bf \Phi}$ are direct modes and ${\bf \Psi}$ are adjoint modes~\cite{rowley:05pod}.}, and they are ordered according to the variance that they capture in the data $\bf A$; if the columns of $\bf A$ are velocity measurements, then the POD modes are ordered in terms of the kinetic energy that they capture.  This variance/energy content is quantified by the entries of the diagonal matrix $\bf\Sigma$, which are called singular values and appear in descending order.

When the size of each snapshot ($n$), is much larger than the number of snapshots ($q$) collected, $n\gg q$, as in high-dimensional fluid systems, there are at most $q$ non-zero singular values, and it is beneficial to use the {\em method of snapshots}~\cite{Sirovich:1987}.  In this method, we solve the following eigenvalue problem:
\begin{equation*}
{\bf A^*A W} = {\bf W \Sigma}_q^2,
\end{equation*}
where ${\bf \Sigma}_q$ is the $q\times q$ upper-left block of $\bf \Sigma$.  It is then possible to find the first $q$ POD modes corresponding to non-trivial singular values by:
\begin{equation*}
{\bf \Psi}_q = {\bf A W \Sigma}_q^{-1}.
\end{equation*}
The snapshots often exhibit low-dimensional phenomena, so that the majority of variance/energy is contained in a few modes, smaller than the number of snapshots collected.  In this case, the POD basis is typically truncated at a pre-determined cut-off value, such as when the columns contain $99\%$ of the variance, so that only the first $r$ modes are kept.  The SVD acts as a filter, and so often the truncated modes correspond to random fluctuations and disturbances.  If the data in the matrix $\bf A$ is generated by a dynamical system (nonlinear system of ordinary differential equations of order $n$),  It is then possible to substitute the truncated POD expansion for the state $\bf U$ into the governing equation and obtain Galerkin projected dynamics on the $r$ basis modes~\cite{HLBR_turb}.  
Recall that we are assuming that the complex systems under consideration exhibits low-dimensional attractors, thus the Galerkin truncation with only a few modes should provide an accurate prediction of the evolution of the system.
Note that it has also been shown recently that it is possible to obtain a \emph{sketched}-SVD by randomly projecting the data initially and then computing the SVD~\cite{Fowler:2009,Gilbert:2012,Qi:2012}.

\section{Methods - Combining $\ell_1$ and $\ell_2$}\label{sec:methods}
\begin{figure}
\begin{center}
\hspace*{.5cm}\begin{overpic}[width=\textwidth]{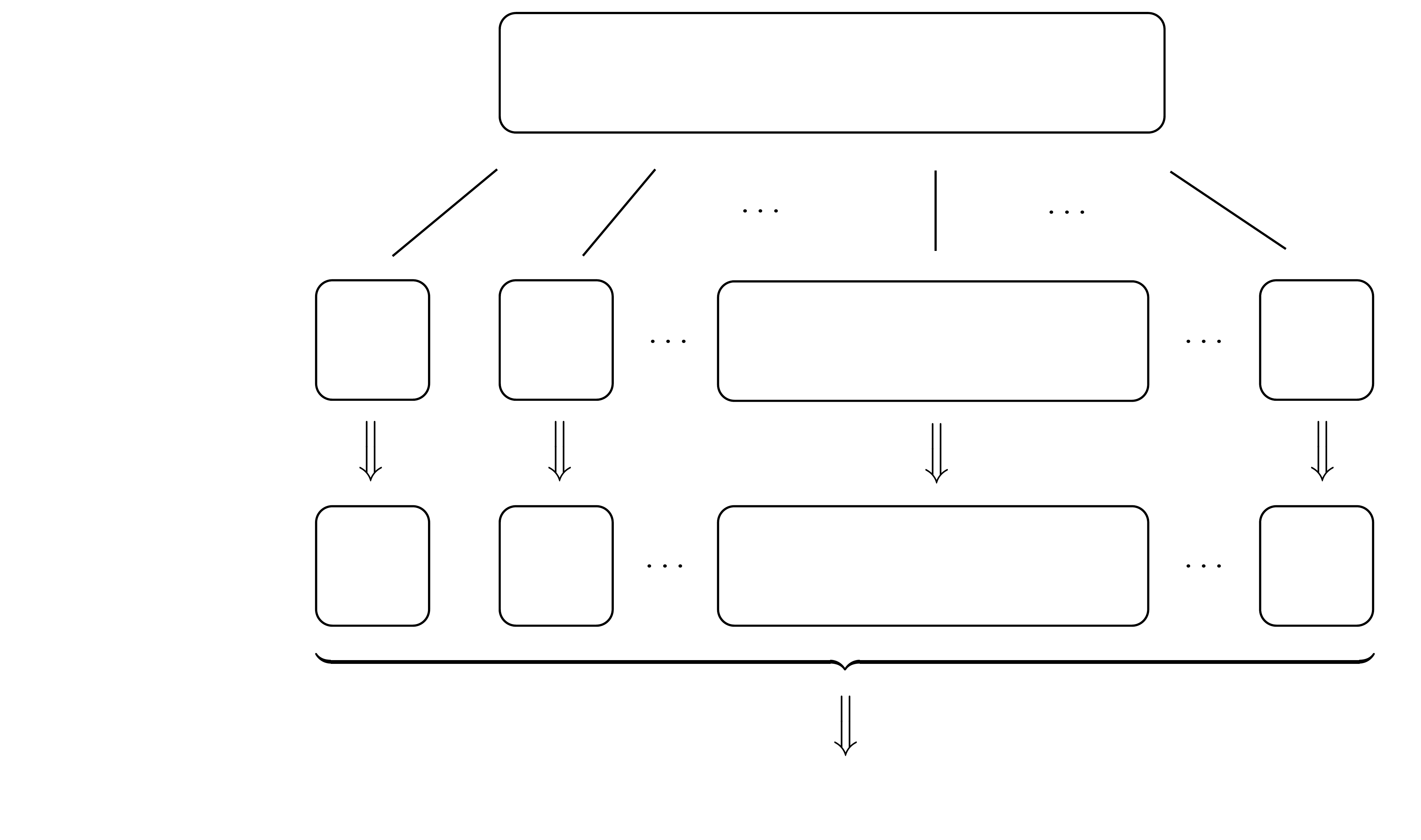} 
\small
\put(46.5,52.5){${\bf U}_t={\bf N}({\bf U},{\bf U}_x,{\bf U}_{xx},\cdots,x,t,\beta)$}
\put(28,45){$\beta_1$}
\put(41.5,45){$\beta_2$}
\put(63.5,45){$\beta_j$}
\put(88.5,45){$\beta_J$}
{\footnotesize
\put(25,34.7){${\bf A}_1$}
\put(38,34.7){${\bf A}_2$}
\put(52.,34.6){{${\bf A}_j =
\begin{bmatrix}
    | &  & |  \\
    {\bf U}(x,t_1)  &\hspace{-.03in} \cdots &\hspace{-.02in} {\bf U}(x,t_q) \\
    | &  & |      
\end{bmatrix}$}}
\put(91.2,34.7){${\bf A}_J$}
\put(25,18.5){${\bf \Psi}_1$}
\put(38,18.5){${\bf \Psi}_2$}
\put(50.8,18.7){${\bf \Psi}_j =
\begin{bmatrix}
    | &  & |  \\
   \psi_1( x,\beta_j) & \hspace{-.06in} \cdots & \hspace{-.05in} \psi_{r_j}(x,\beta_j) \\ 
    | &  & |      
\end{bmatrix}$}
\put(91.2,18.5){${\bf \Psi}_J$}
}
\put(49.5,3){${\bf \Psi}=\begin{bmatrix}{\bf \Psi}_1 & {\bf \Psi}_2 & \cdots & {\bf \Psi}_J \end{bmatrix}$}
{\bf\put(36.5,56){Complex System}}
{\bf\put(-3,46){Step 1:}
\put(-3,44){Collect Data}}
\put(-3,42){(simulations or experiments)}
{\normalfont\bf\put(-3,28){Step 2: }
\put(-3,26){Dimension Reduction}
\put(-3,9){Step 3:}
\put(-3,7){Library Building}}
\end{overpic}
\caption{Schematic of $\ell_2$ strategy for library building.  Data ${\bf A}_j$ are collected for many values of the bifurcation parameter $\beta_j$, and the principal components of this data are computed and truncated in ${\bf\Psi}_j$.  Although each basis ${\bf \Psi}_j$ is truncated to contain only the most energetically relevant structures, the concatenated library ${\bf \Psi}$ is overcomplete.}
\label{fig:L2schematic}
\end{center}
\end{figure}

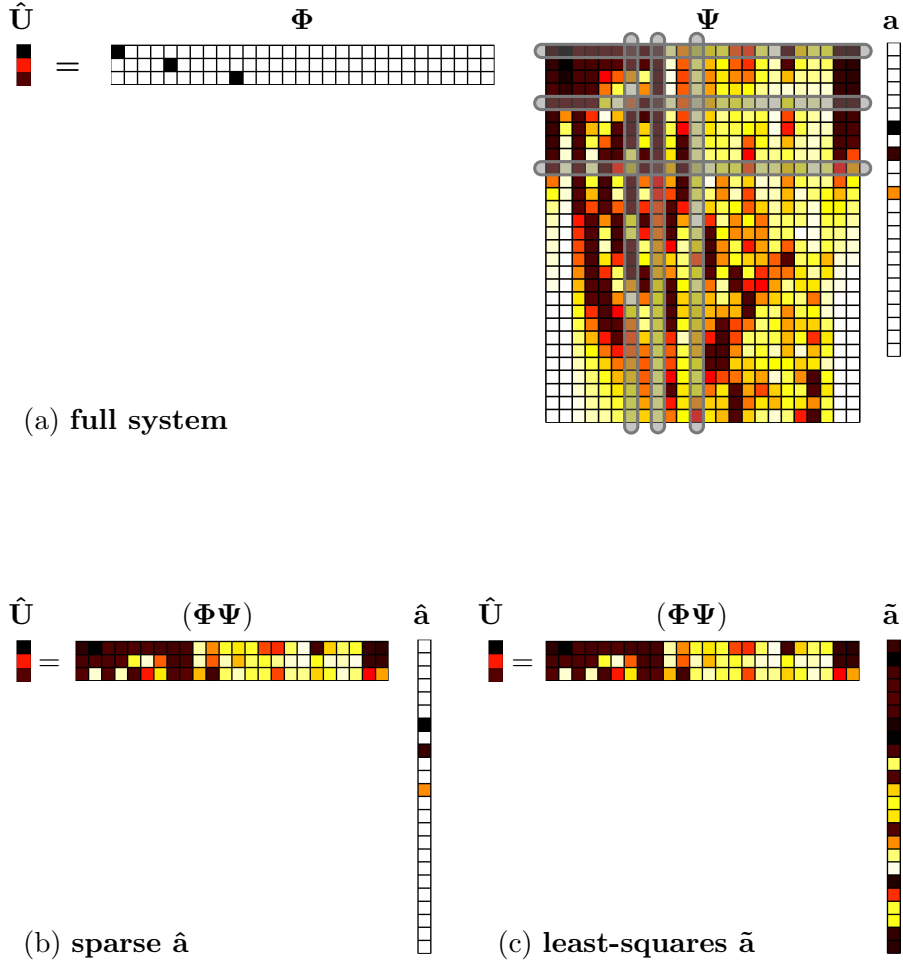
\begin{figure}
\begin{center}
\begin{overpic}[width=.75\textwidth]{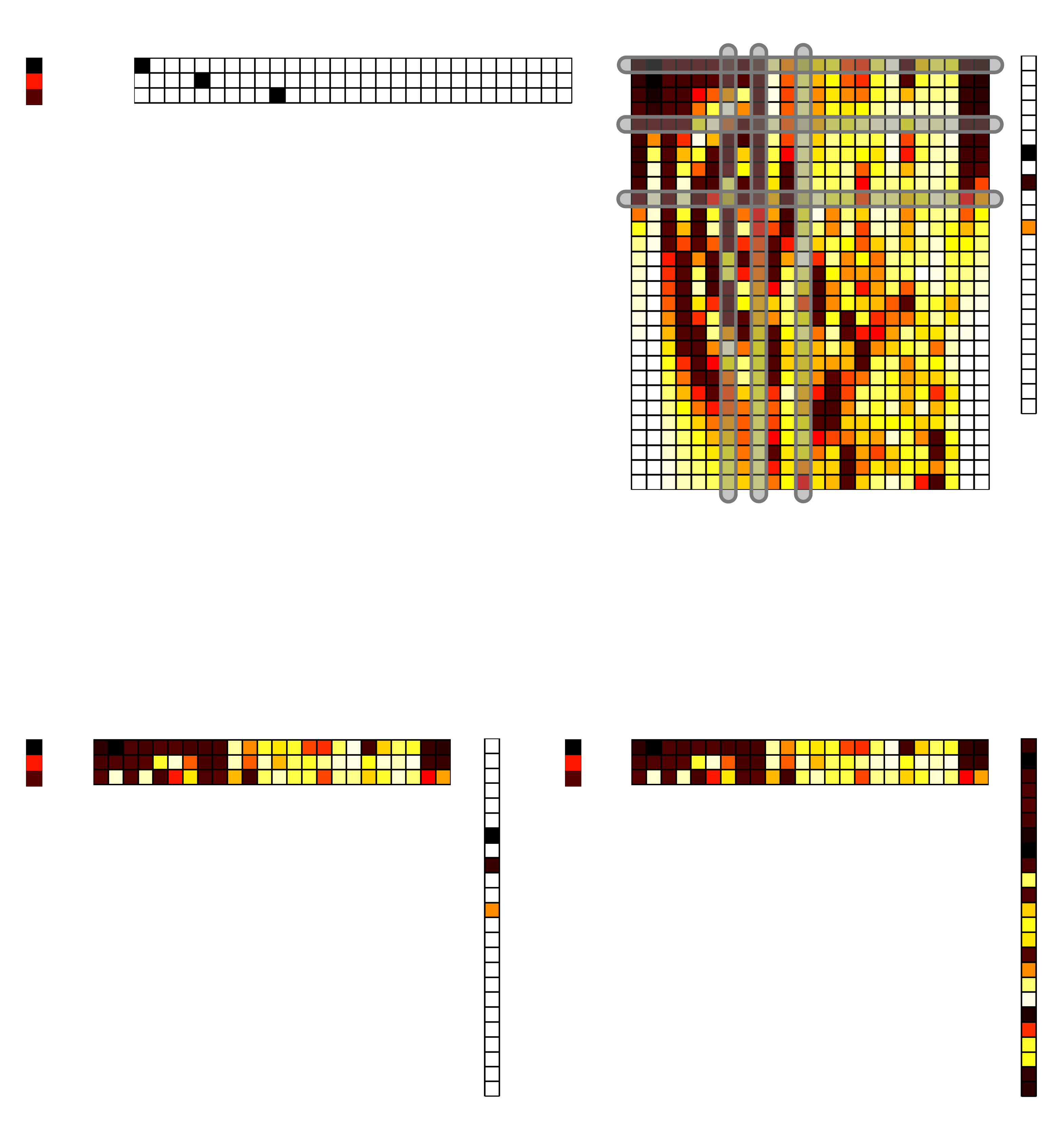}
\normalsize
\put(1.5,96.5){$\bf \hat U$}
\put(30,96.5){${\bf \Phi}$}
\put(71,96.5){${\bf \Psi}$}
\put(90,96.5){${\bf a}$}
\put(6.2,91.9){$\boldsymbol{=}$}
\put(3,56){(a) {\bf full system}}
\put(1.5,36.25){$\bf \hat U$}
\put(19,36.25){$({\bf \Phi\Psi})$}
\put(42.5,36.25){${\bf \hat a}$}
\put(4.5,31.6){$=$}
\put(3,3){(b) {\bf sparse ${\bf \hat a}$}}
\put(49.,36.25){$\bf \hat U$}
\put(67,36.25){$({\bf \Phi\Psi})$}
\put(90,36.25){${\bf \tilde a}$}
\put(52.4,31.6){$=$}
\put(51,3){(c) {\bf least-squares ${\bf \tilde a}$}}
\end{overpic}
\end{center}
\caption{Schematic of identification of sparse mode amplitudes $\bf\hat a$ by $\ell_1$-minimization.  (a) Illustration of measurement matrix $\bf \Phi$ and sparse basis $\bf \Psi$.  The underdetermined matrix $({\bf \Phi\Psi})$ admits a sparse solution ${\bf \hat a}$ (b) and a least-squares solution ${\bf \tilde a}=({\bf\Psi^*\Phi^*\Phi\Psi})^{-1}{\bf \Psi^*\Phi^*\hat U}$ (c).  This type of diagram was introduced by Baraniuk in~\cite{Baraniuk:2007}.  The data used in this figure is from the CQGLE system in Sec.~\ref{sec:results}.}
\label{fig:L1schematic}
\end{figure}

The major contribution of this work is the combination of library building techniques (depicted schematically in Figure~\ref{fig:L2schematic}) based on the $\ell_2$-optimal proper orthogonal decomposition (POD) with the $\ell_1$-based compressive sensing (CS) architecture (depicted schematically in Figure~\ref{fig:L1schematic}) for classification and reconstruction.   

Consider a complex system that evolves according to the partial differential equation
\begin{equation}
  {\bf U}_t = {\bf N} ({\bf U},{\bf U}_x,{\bf U}_{xx}, \cdots , x, t, \beta),
  \label{eq:U_beta}
\end{equation}
where ${\bf U}(x,t)$ is a vector of physically relevant quantities and the subscripts $t$ and $x$ denote partial differentiation in time and space, respectively.  Note that higher-spatial dimensions may be considered without loss of generality.
The function ${\bf N}(\cdot)$ can be a complicated, nonlinear function of the quantity ${\bf U}$, its derivatives, and both space and time.
The parameter $\beta$ is a bifurcation parameter with respect to which the solution of the governing PDE changes markedly.
We assume a spatial discretization of (\ref{eq:U_beta}), which yields a high-dimensional system of degree $n$.

We would like to use measurements of the system (\ref{eq:U_beta}) to determine its state.
However, full-state measurements are impractical for the high-dimensional system generated by discretization.
Instead, $m$ measurements are taken, where $m\ll n$; thus the measurements are {\em sparse}.
In this work, we consider {\em spatially localized} or \emph{point} measurements, $\bf \hat U$, as discussed in Sec.~\ref{sec:back:l1}.  In this case, the matrix ${\bf\Phi}\in\mathbb{R}^{m\times n}$ from \eqref{eq:sparse2} is comprised of rows of the identity matrix corresponding to the measurement locations.  These $m$-dimensional sparse observations are used to reconstruct the full $n$-dimensional state vector ${\bf U}$.  

Our approach is to {\em learn} a library of low-rank dynamical approximations in which the dynamics are {\em sparse} and then apply compressive sensing to reconstruct the dynamics from $m\ll n$ measurements.  
First, we explore the full system (\ref{eq:U_beta}) and collect dense measurements for various values $\beta_1,\beta_2,\cdots, \beta_J$ of interest, making sure to cover a number of unique dynamical regimes.  
For each case, snapshots of data from simulations or experiments are taken at a number of instances in time and organized into a data matrix describing the evolution of the full-state system: 
\begin{equation*}
{\bf A}_j =
\begin{bmatrix}
    | & | & & |  \\
    {\bf U}(x,t_1)  & {\bf U}(x,t_2) & \cdots & {\bf U}(x,t_q) \\
    | & | & & |      
\end{bmatrix},
\end{equation*}
where $q$ is the number of snapshots taken.

Once the data matrix is constructed for a given $\beta_j$, its POD modes, or principal components, ${\bf \Psi}_j$ are identified through a singular value decomposition (SVD):  \mbox{${\bf \Psi}_j=\left\{\psi_i(x,\beta_j)\right\}_{i=1}^{r_j}$}.
The POD modes are orthogonal and ordered by energy content.
The number of modes retained, $r_j$, is determined by a cut-off criterion; for instance, one might specify that modes comprising 99\% of the energy be kept.

With the modes identified for each $\beta_j$, an overcomplete library ${\bf \Psi}$ is constructed that contains all of the low-rank approximations for each dynamic regime:
\begin{equation}
  \hspace*{-.1in}{\bf\Psi}\!\!=\!\! \left[ \!\!\! \begin{array}{ccccccccc}  | &   & | &  & | &  & | & \\
        \psi_1\!(x,\!\beta_1\!)\!\! &  \!\! \cdots \!\!& \!\! \psi_{r_1}\!(x,\!\beta_1\!)\!\! &
       \!\!\cdots \!\!&
       \!\! \psi_1\!(x,\!\beta_J\!)\!\! &  \!\! \cdots \!\!& \!\! \psi_{r_J}\!(x,\!\beta_J\!)\!\!  \\
        | &   & | &  & |&  &| &  \end{array}   \! \!\!\!\!\! \right].
        \label{eq:library}
\end{equation}
The library ${\bf\Psi}\in\mathbb{R}^{n\times p}$ contains the representative low-rank modes for all of the dynamical behavior of the governing system that we explored in simulations or experiments.
This is the {\em supervised learning} portion of the analysis, resulting in a small number ($p\ll n$) of library elements; note that $p=\sum_{j=1}^Jr_j$.  
The $p$ library modes are {\em not} orthogonal, but rather come in groups of orthogonal POD modes for each given $\beta_j$.  
The dynamics at any given time will belong to a specific $\beta_j$ regime so that that the instantaneous dynamics are {\em sparse} in the library basis, allowing for a sparse representation~\cite{Wright:2009}.  
This \emph{overcomplete library building} procedure is summarized in Figure~\ref{fig:L2schematic}.

With the library (\ref{eq:library}), we can expand the state $\bf{U}$ using the low-rank POD representation
\begin{equation}
  {\bf U}(x,t)=\sum_{j=1}^{J} \sum_{r=1}^{r_j} a_{jr} (t) \psi_r (x,\beta_j) = {\bf\Psi} {\bf a}.
\label{eq:fullPOD}
\end{equation}
The solution is now represented in the $p$ library elements constructed for the various values of $\beta$, and by construction, we expect $\bf{a}$ to be sparse in the basis ${\bf \Psi}$.  
This is because for any particular $\beta_j$, only a small subset of library elements is required to represent the solution.

Equation~\eqref{eq:fullPOD} is of the form in \eqref{eq:sparse1}.  
To determine the vector $\bf a$ from a sparse data measurement ${\bf \hat{U}}={\bf \Phi}\bf{U}$, insert \eqref{eq:fullPOD} into \eqref{eq:sparse2} and solve the under-determined linear system ${\bf \hat{U}} = ({\bf\Phi \Psi}) {\bf a}$ from \eqref{eq:sparse3} which has $m$ equations (constraints) and $p$ unknowns (modal coefficients), with  $m\ll p$.  
We solve for a sparse $\bf{\hat a}$ using compressive sensing ($\ell_1$ minimization).  
This approach is natural because it promotes sparsity, an expected property of $\bf{a}$.  
Further, solving for $\bf{a}$ using $\ell_1$ minimization in the reduced-order library basis is significantly more efficient than solving for $\bf{U}$ in the full space since $p\ll n$.  
The sparsity-promoting compressed sensing procedure is illustrated in Figure~\ref{fig:L1schematic}.

The library construction (a one-time cost)  and sparse sensing combine to give an efficient algorithm for approximating the low-rank dynamics of the full PDE (\ref{eq:U_beta}) using a limited number of sensors and an empirically determined, overcomplete database.  
Specifically, the full-state of the system $\bf{U}$ at any given time $t$ is achieved by evaluating $\bf{a}$.  
There are a number of immediate advantages to this method for characterizing complex dynamical systems:
\begin{itemize}
\item[(i)] Once the library is constructed from extensive simulations, future prediction of the system is efficient since the correct POD modes for any dynamical regime $\beta_j$ have already been computed,  
\item[(ii)]  The algorithm works equally well with experimental data
 in an equation-free context, for instance by using dynamic mode decomposition \cite{Rowley:2009,schmid:2010} or equation-free modeling \cite{Kevrekidis:2003} in place of POD,   
\item[(iii)] Given the low-rank space in which the algorithm works, it is ideal for use with control strategies, which are only practical for real-time application with low-dimensional systems.
\end{itemize}

Sparse sensing is significantly less expensive in the learned library $\bf\Psi$ since the high-dimensional state has been replaced with a truncated POD representation.  
Additionally, less information is required to categorize a signal than is required to fully reconstruct the signal, as in the compressive sensing paradigm.  
This combination of classification and reconstruction in a concatenated set of truncated POD bases using $\ell_1$ minimization is appealing on a number of levels.  
There is also benefit to keeping the individual POD bases ${\bf\Psi}_j$ for reconstruction once the bifurcation regime $\beta$ has been identified.

\section{Results}\label{sec:results}
To illustrate the aforementioned strategy, consider the complex Ginzburg-Landau model \cite{Cross:1993}, which is ubiquitous in mathematical physics.  
Here it is modified to include both quintic terms and a fourth-order diffusion term much like the Swift-Hohenberg equation:
\begin{eqnarray}
 i {\bf U}_t \!+\! \left( \! \frac{1}{2} \!-\! i\tau \!\! \right) \!{\bf U}_{xx} \!-\! i \kappa {\bf U}_{xxxx}
 + (1\!-\! i\mu) |{\bf U}|^2 {\bf U} + (\nu \!-\! i \varepsilon) |{\bf U}|^4 {\bf U} \!\!-\! i\gamma {\bf U} \! =\! 0,
  \label{eq:cqgle_ml}
\end{eqnarray}
where ${\bf U}(x,t)$ is a complex function of space and time.  
Interesting solutions to this governing equation abound, characterized by the parameter values
$ {\beta}=(\tau, \kappa, \mu, \nu, \varepsilon, \gamma)$.  
In particular, we consider six regimes that illustrate different dynamical behaviors, described in Table~\ref{ta:cqgle_beta}.

\begin{table}[t]
\begin{center}
\caption{\label{ta:cqgle_beta} Parameter regimes $\beta_j$ for the complex Ginzburg-Landau equation (\ref{eq:cqgle_ml})
(See Figure~\ref{fig:cqgle1}).  The low-rank approximations of these parameter regimes are used to construct the elements of
the library ${\bf \Psi}$.}
  \begin{tabular}{|c|c|c|c|c|c|c|c|}
  \hline
         &  $\tau$ & $\kappa$ & $\mu$ & $\nu$ & $\varepsilon$ & $\gamma$  & description \\ \hline
     $\beta_1$       &  -0.3 & -0.05 & 1.45 & 0 & -0.1 & -0.5  & 3-hump, localized  \\ \hline
      $\beta_2$       &  -0.3 & -0.05 & 1.4 & 0 & -0.1 & -0.5 & localized, side lobes\\ \hline
      $\beta_3$        &  0.08 & 0 & 0.66 & -0.1 & -0.1 & -0.1 & breather  \\ \hline
       $\beta_4$      &  0.125 & 0 & 1 & -0.6 & -0.1 & -0.1 & exploding soliton  \\ \hline
     $\beta_5$        &  0.08 & -0.05 & 0.6 & -0.1 & -0.1 & -0.1 & fat soliton  \\ \hline
    $\beta_6$          &  0.08 & -0.05 & 0.5 & -0.1 & -0.1 & -0.1 & dissipative soliton  \\\hline
  \end{tabular}
\end{center}
\end{table}

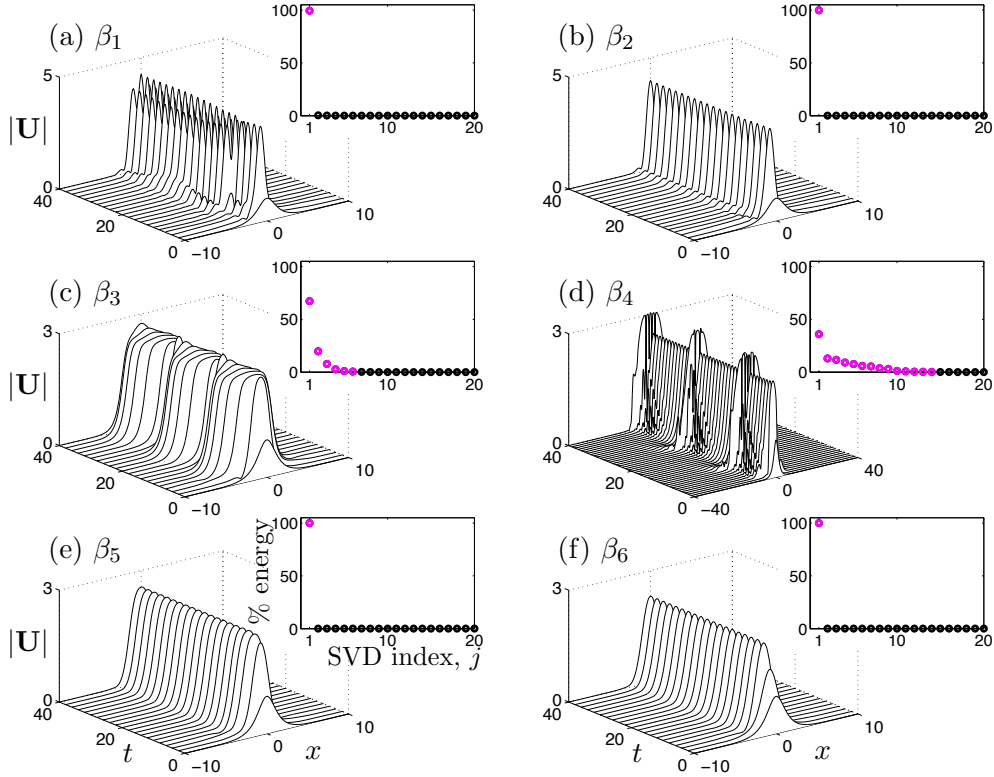
\begin{figure}
\begin{center}
\begin{tabular}{lr}
\hskip -.65in
\begin{overpic}[width=.3\textwidth]{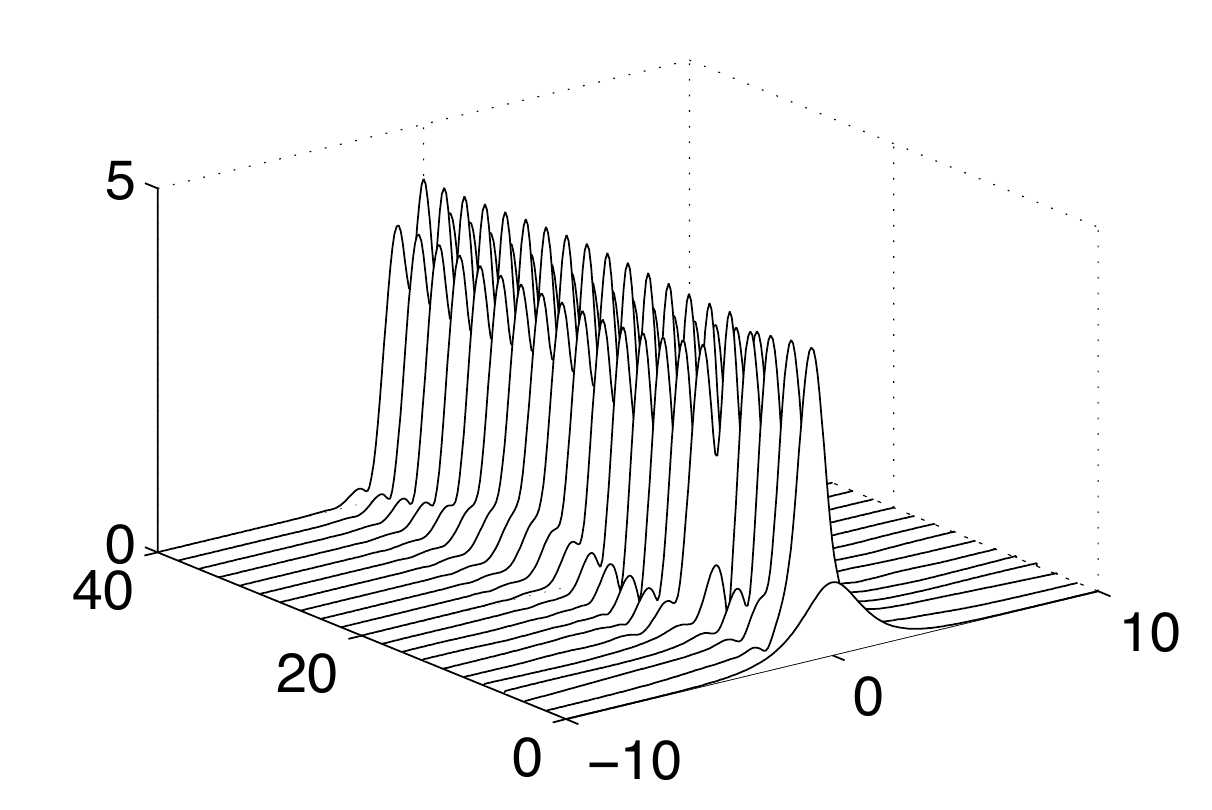} 
\put(70,35){\includegraphics[width=.18\textwidth]{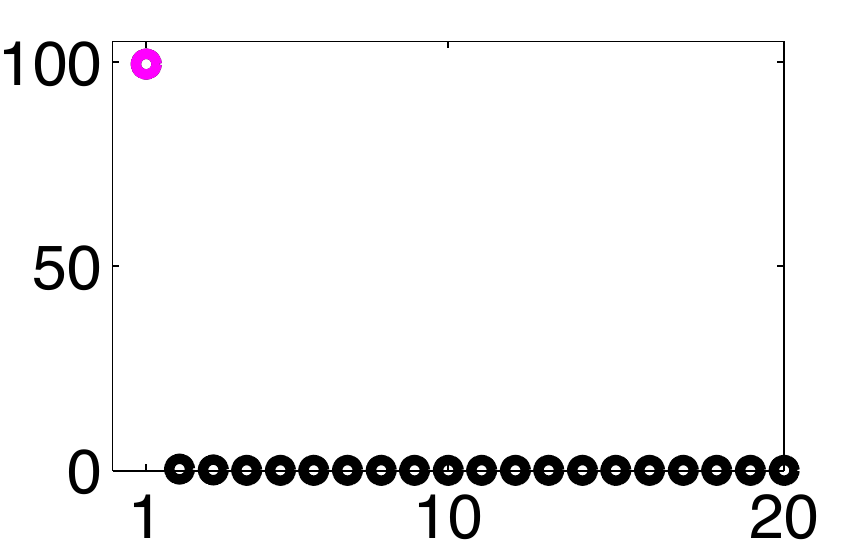}}
\normalsize
\put(10,60){(a) $\beta_1$}
\put(-1,35){$|{\bf U}|$}
\end{overpic}&\hskip .55in
\begin{overpic}[width=.3\textwidth]{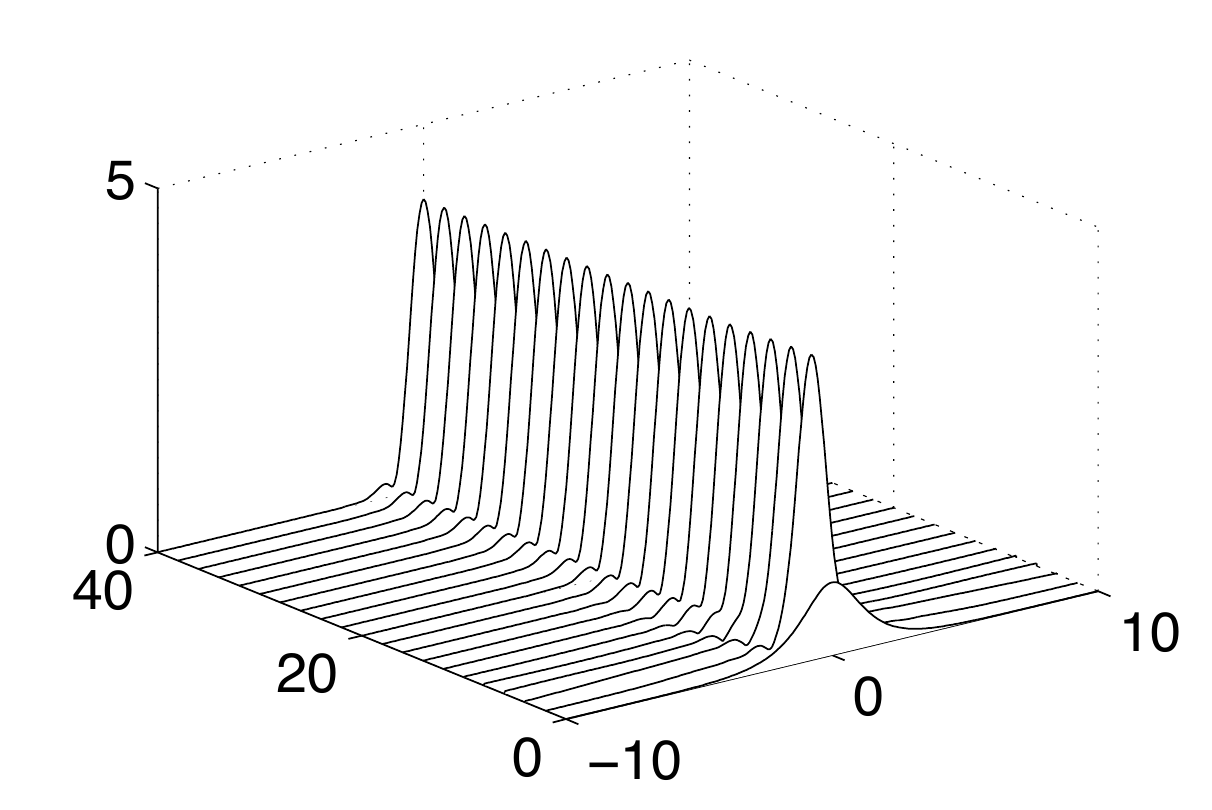} 
\put(70,35){\includegraphics[width=.18\textwidth]{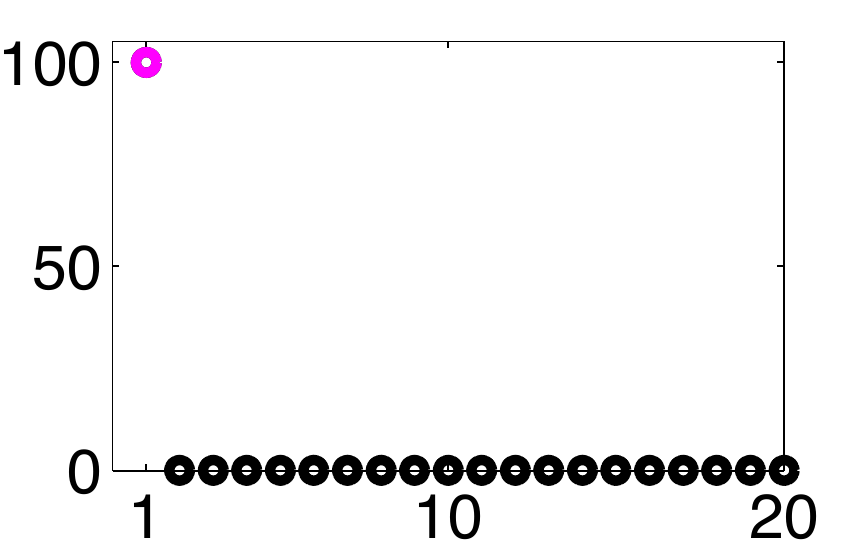}}
\normalsize
\put(10,60){(b) $\beta_2$}
\end{overpic}\\
\hskip -.65in
\begin{overpic}[width=.3\textwidth]{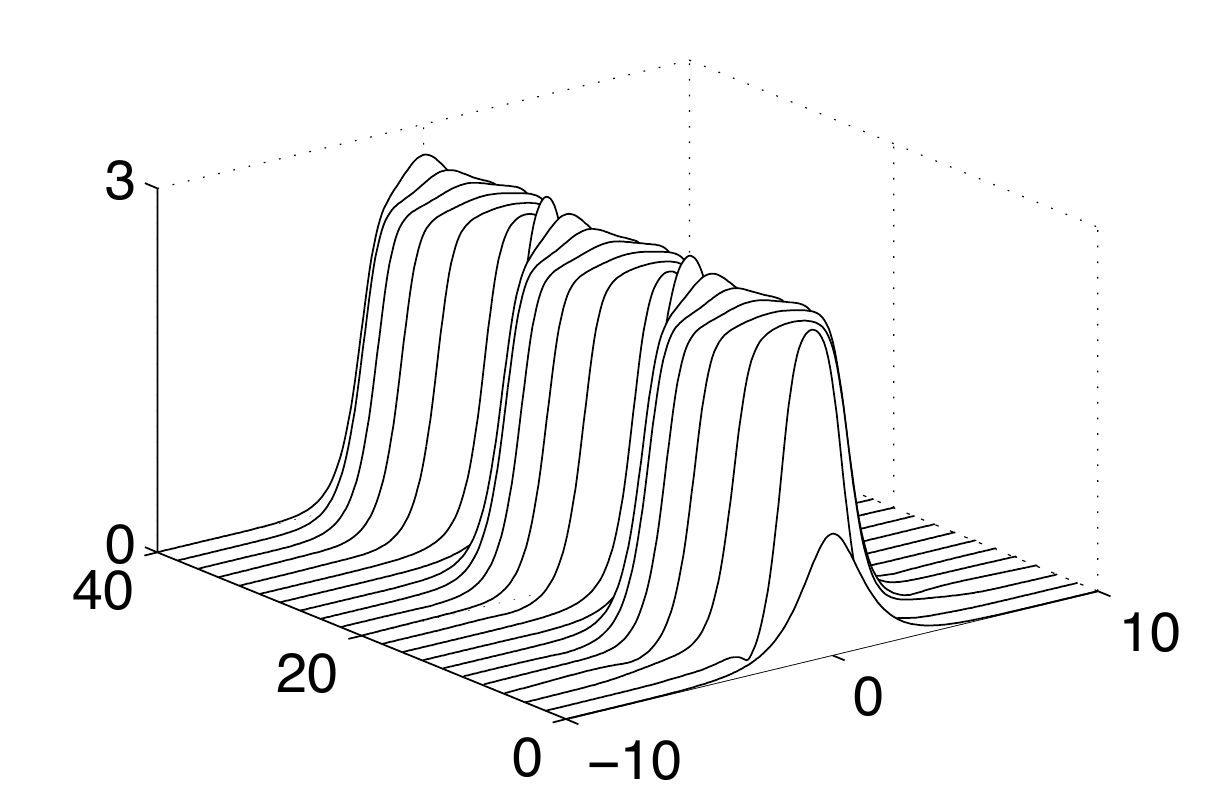}  
\put(70,35){\includegraphics[width=.18\textwidth]{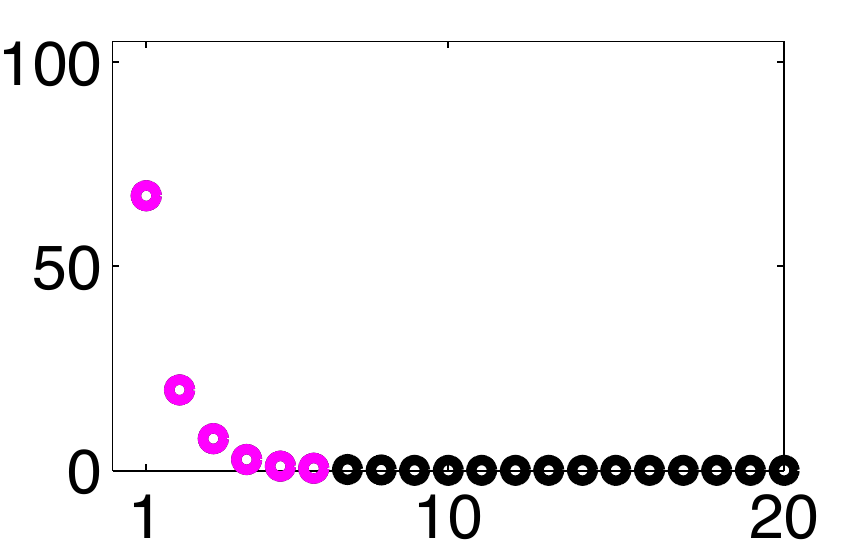}}
\normalsize
\put(10,60){(c) $\beta_3$}
\put(-1,35){$|{\bf U}|$}
\end{overpic}&\hskip .55in
\begin{overpic}[width=.3\textwidth]{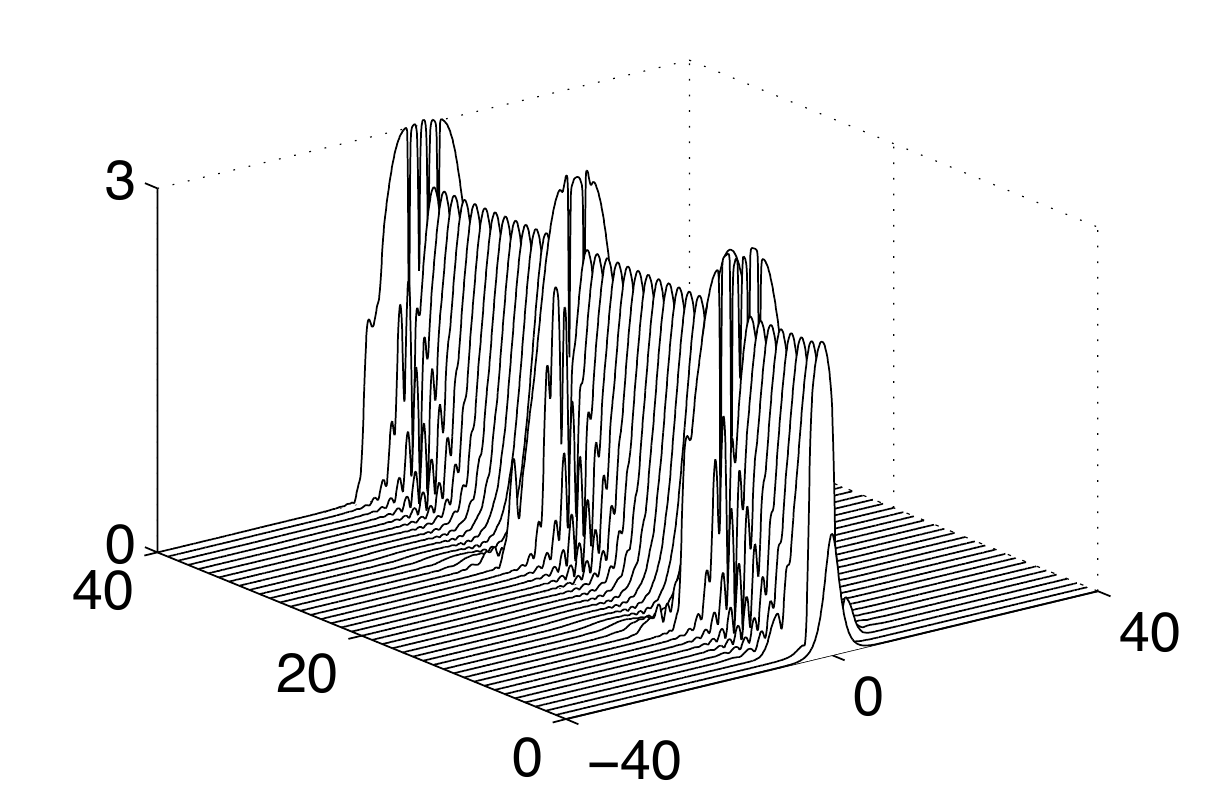} 
\put(70,35){\includegraphics[width=.18\textwidth]{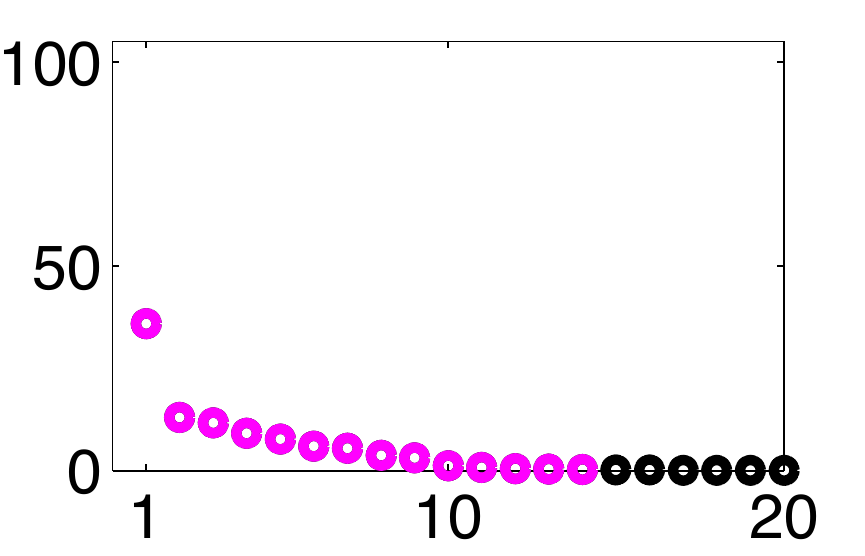}}
\normalsize
\put(10,60){(d) $\beta_4$}
\end{overpic}\\  
\hskip -.65in
\begin{overpic}[width=.3\textwidth]{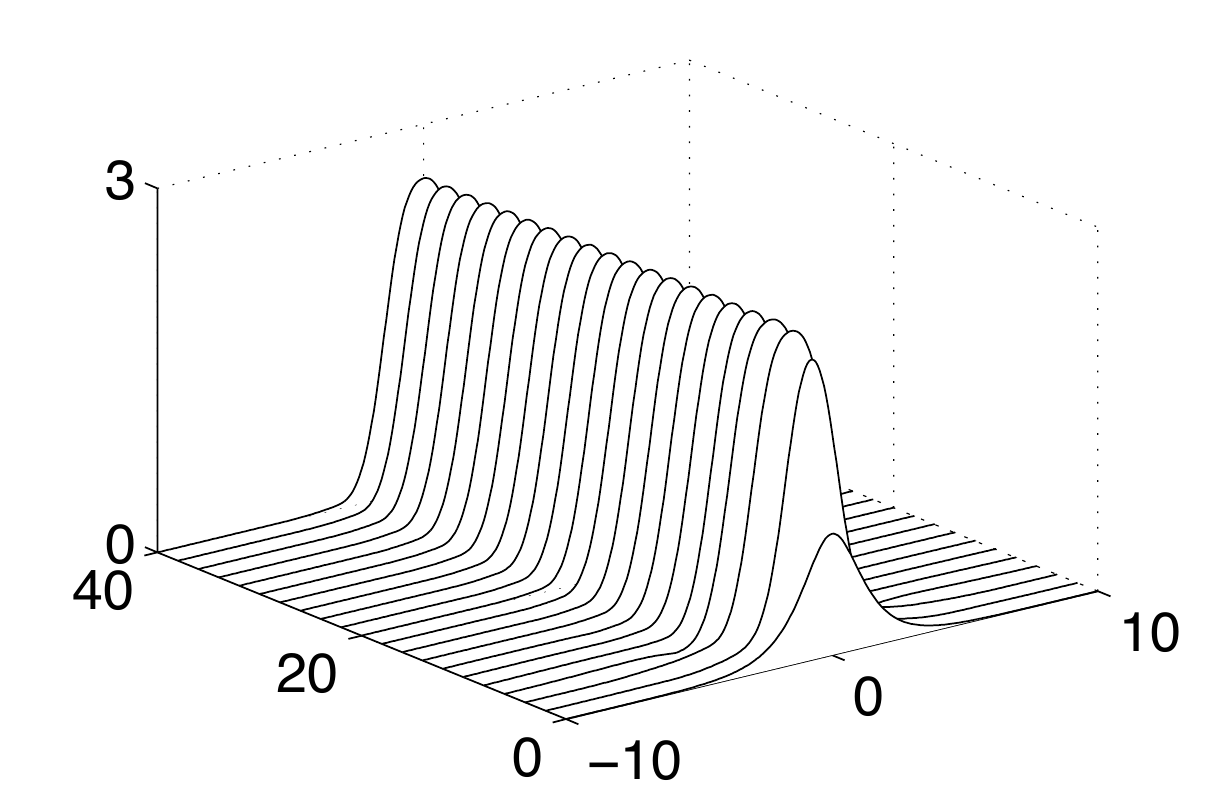}  
\put(70,35){\includegraphics[width=.18\textwidth]{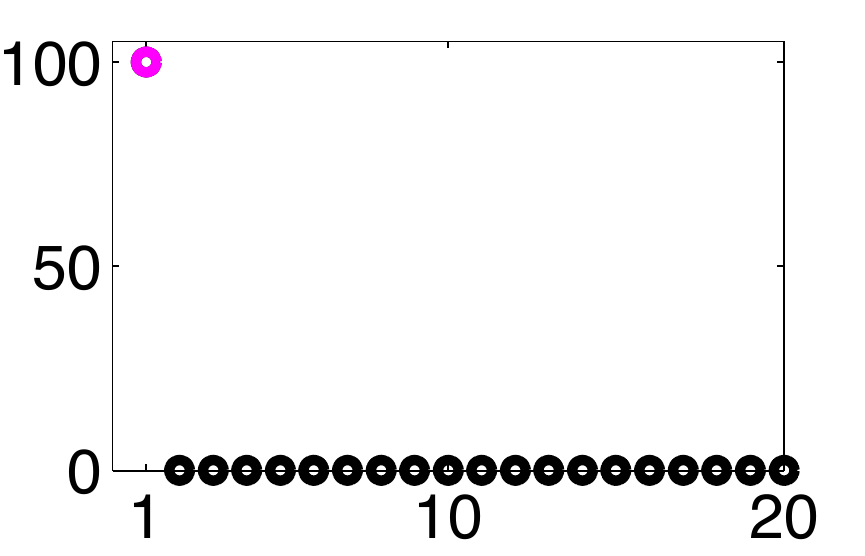}}
\normalsize
\put(10,60){(e) $\beta_5$}
\put(-1,35){$|{\bf U}|$}
\put(30,5){$t$}
\put(79,6){$x$}
\put(85,31){\small SVD index, $j$}
\put(69,42){\small\begin{rotate}{90}\% energy\end{rotate}}
\end{overpic}&\hskip .55in
\begin{overpic}[width=.3\textwidth]{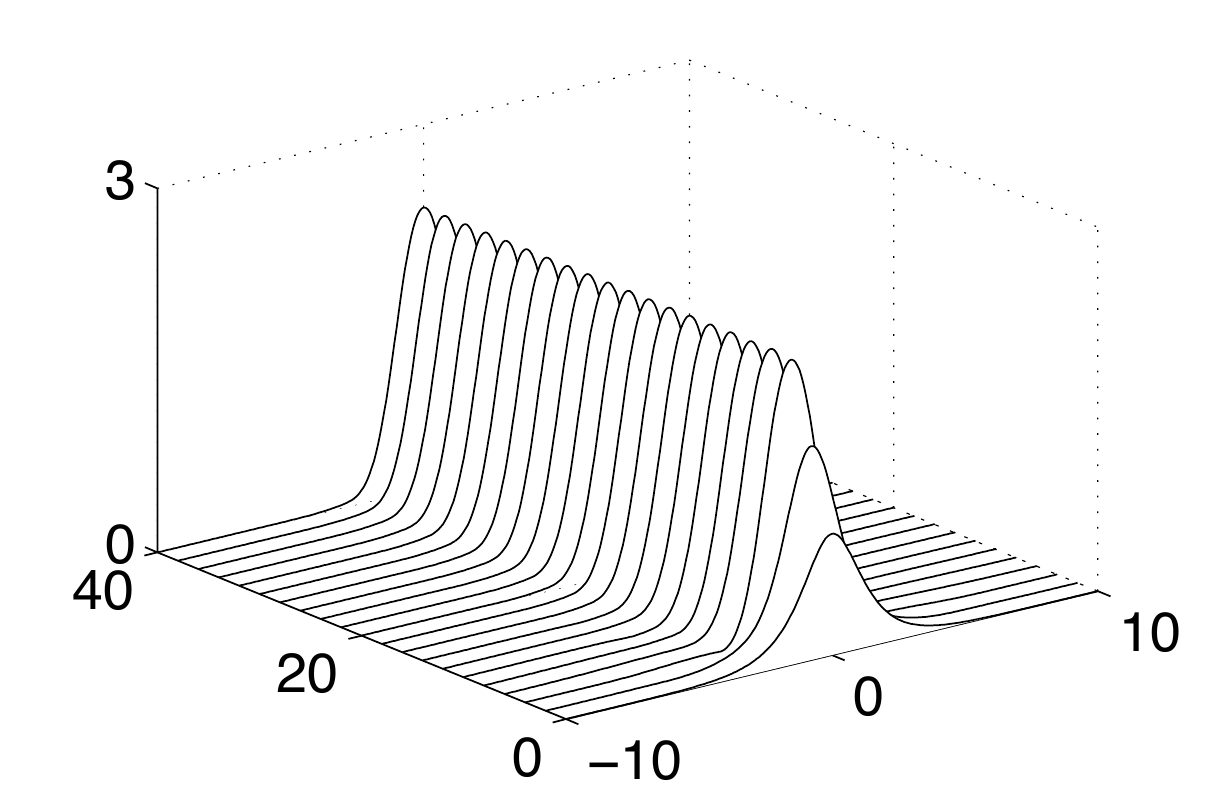} 
\put(70,35){\includegraphics[width=.18\textwidth]{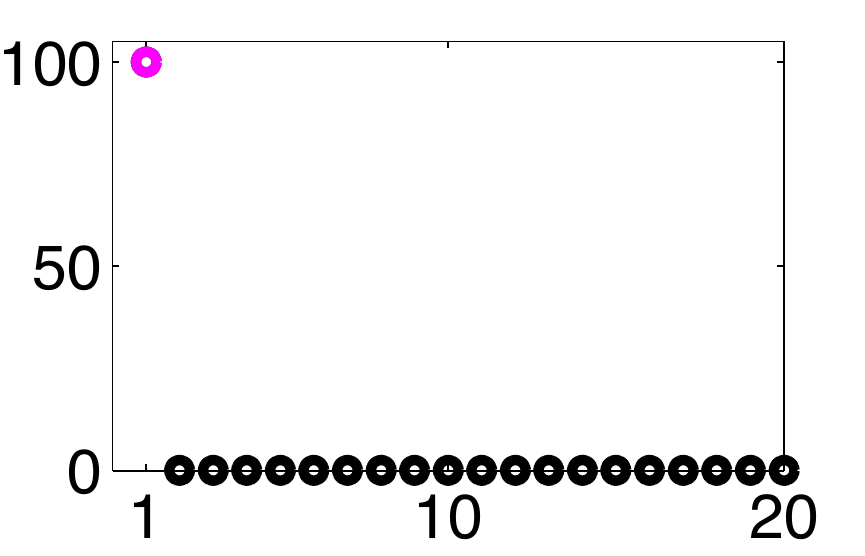}}
\normalsize
\put(10,60){(f) $\beta_6$}
\put(30,5){$t$}
\put(79,6){$x$}
\end{overpic}\\
\end{tabular}
\end{center}
\caption{\label{fig:cqgle1} Evolution dynamics of (\ref{eq:cqgle_ml}) for the six parameter
regimes given in Table~\ref{ta:cqgle_beta}: (a) $\beta_1$, (b) $\beta_2$, (c) $\beta_3$, (d) $\beta_4$,
(e) $\beta_5$, and (f) $\beta_6$.  All parameter regimes exhibit stable, low-dimensional attractors as evidenced by the singular values (inset).  The SVD sampling occurs for every $\Delta t=1$ in the interval $t\in[40,80]$.  Magenta circles represent the modes that comprise 99\% of the data and are used for the library ${\bf \Psi}$.}
\end{figure}

Figure~\ref{fig:cqgle1} illustrates the corresponding low-rank behavior produced in the simulations.  

\begin{figure}[t]
\begin{center}
\begin{overpic}[width=.9\textwidth]{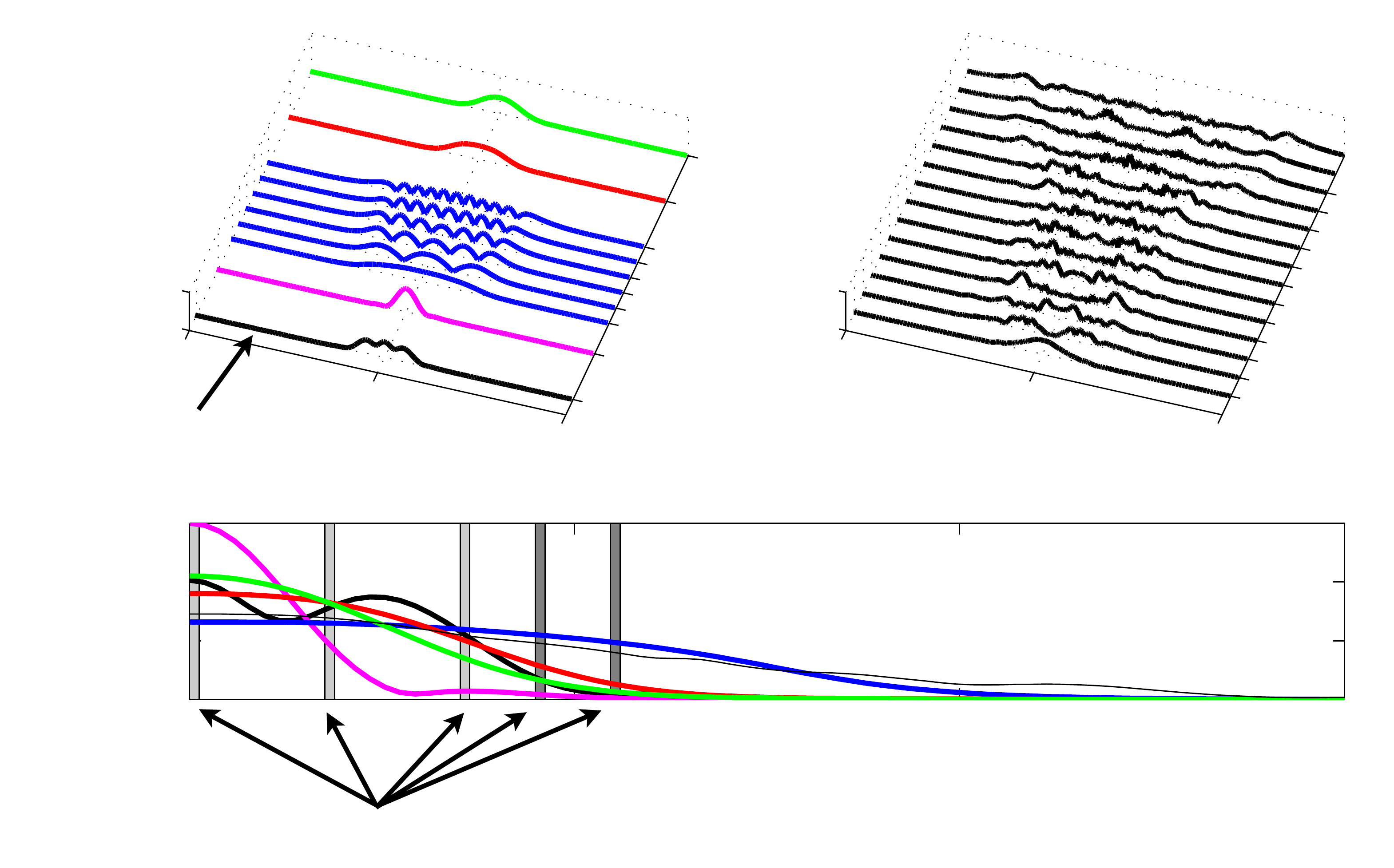}
\normalsize
\put(4,38){$|{\bf U}|$}
\put(9,40.6){0.5}
\put(11,38){0}
\put(9,35){-10}
\put(25,32){0}
\put(37.5,29){10}
\put(31,28){$x$}
\put(35,58){\bf library modes}
\put(10,30.5){$\psi_1(x,\beta_1)$}
\put(42,32){(a) $\beta_1$}
\put(43.8,35.5){(b) $\beta_2$}
\put(9.6,48.5){(c) $\beta_3$}
\put(12.,53.5){(e) $\beta_5$}
\put(14,57.5){(f) $\beta_6$}
\put(52.7,38){$|{\bf U}|$}
\put(56,40.6){0.4}
\put(58,38){0}
\put(57,35){-20}
\put(72,32){0}
\put(84.5,29){20}
\put(78,28){$x$}
\put(93,40){(d) $\beta_4$}
\put(11.5,11.5){0}
\put(10,15.5){0.1}
\put(10,19.5){0.2}
\put(10,23.5){0.3}
\put(4,20){$|{\bf U}|$}
\put(13,8.5){0}
\put(40.6,8.5){2}
\put(68,8.5){4}
\put(95.8,8.5){6}
\put(83,7){$x$}
\put(28.5,2.5){\bf sensor locations}
\end{overpic}
\end{center}
\caption{\label{fig:cqgle2} Library ${\bf \Psi}$ of the dominant modes.  The groupings, identified by their $\beta_j$ value, are associated with the different dynamical regimes (a)-(f) in Figure~\ref{fig:cqgle1}.  Note that the modes of the exploding dissipative solution (d) have been including separately in the right panel as there are 14 modes
  required to capture 99\% of the dynamics in this regime.  A sample cross-section
  of the first mode of each library element $\psi_1(x,\beta_j)$ (j=1,2,3,4,5,6) is shown
  in the bottom panel color-coded with the top panels. 
  The bottom panel also shows the spatial location of the three sensors (light
  shaded) and five sensors (addition dark shaded) used for sparse sampling. }
\end{figure}

As is common in many complex dynamical systems, especially those of a dissipative nature, low-dimensional
attractors are embedded in the high-dimensional space.  
The simulations from each of these dynamic regimes exhibit low-dimensional structures which are spontaneously formed from generic, localized initial data.  
The low-dimensional structures allow for the low-rank POD approximations used in the library construction of Figure~\ref{fig:cqgle2}, as described in \eqref{eq:library} and Figure~\ref{fig:L2schematic}.

\begin{figure}[t]
\begin{center}
\begin{overpic}[width=.95\textwidth]{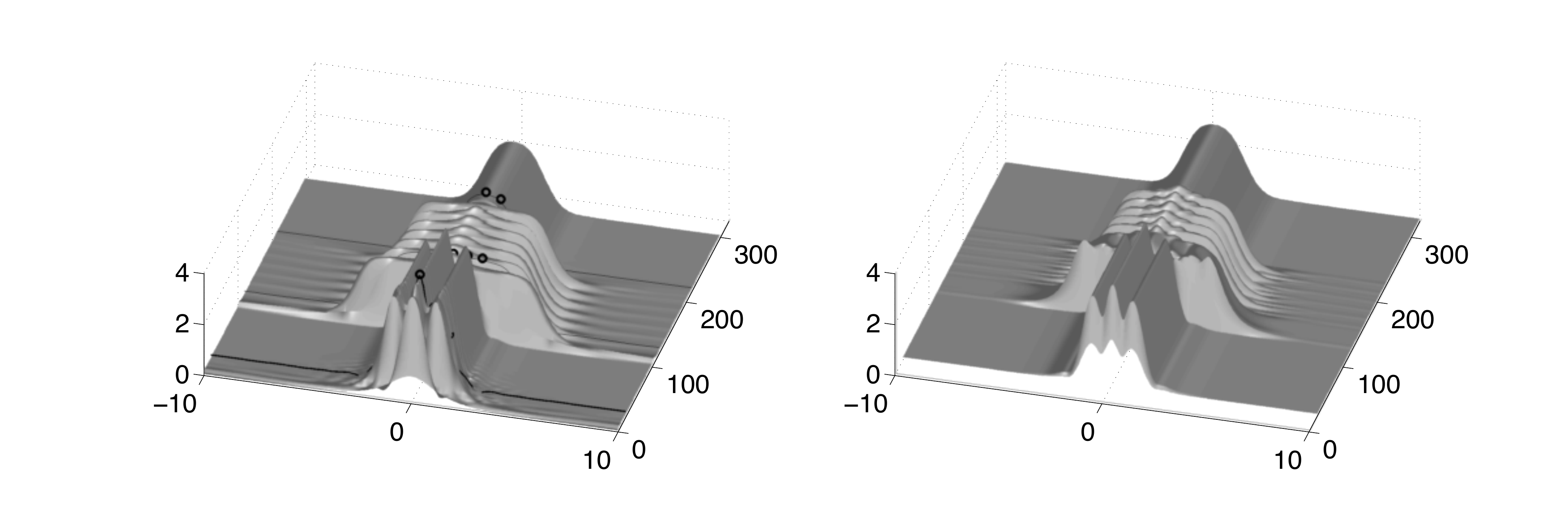}
\normalsize
\put(7,30){(a)}
\put(53,30){(b)}
\put(12,30){\bf Exact simulation}
\put(58,30){\bf Dynamic reconstruction}
\put(30,2){$x$}
\put(74,2){$x$}
\put(7,11){$|{\bf U}|$}
\put(51,11){$|{\bf U}|$}
\put(47,8.5){$t$}
\put(91,8.5){$t$}
\end{overpic}
\\
\vskip .2in
\begin{overpic}[width=.95\textwidth]{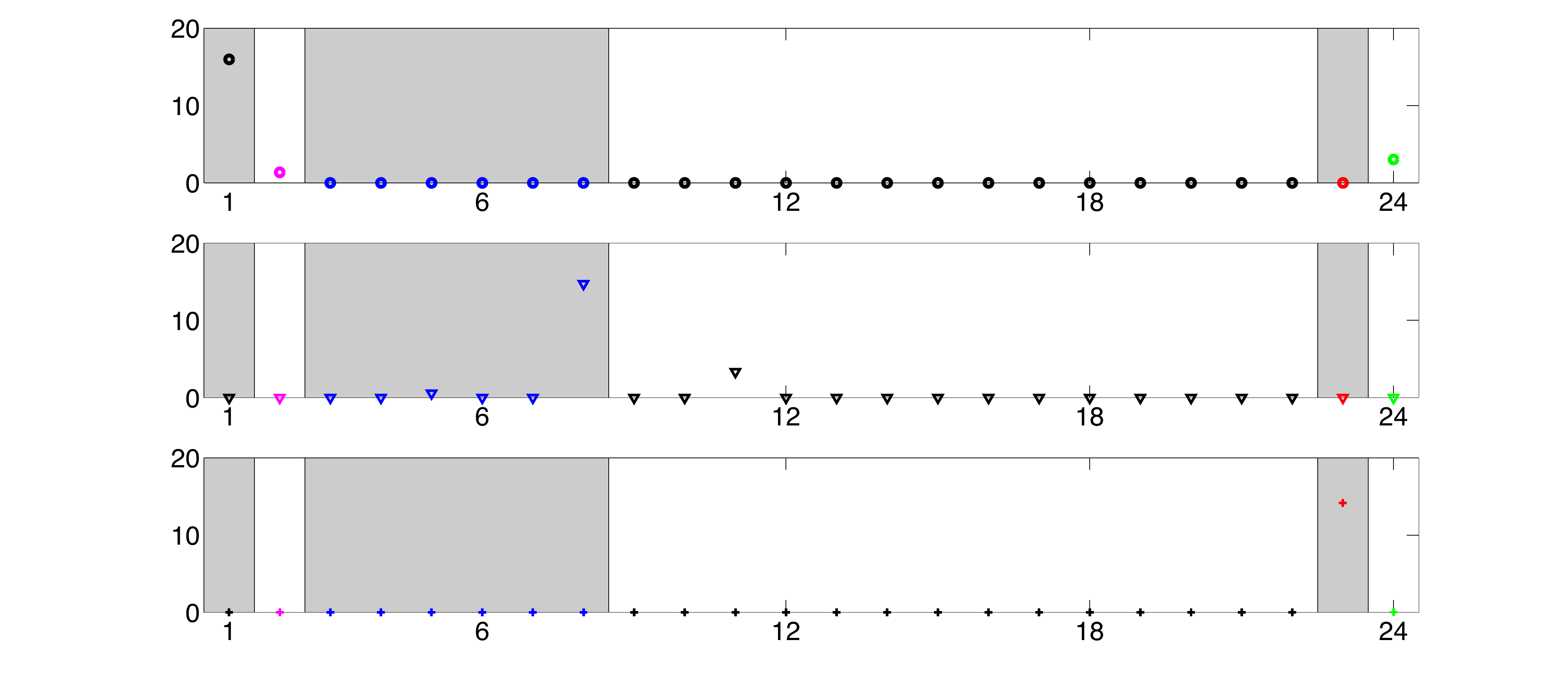}
\normalsize
\put(7,45){(c)}
\put(13.1,44){$\beta_1$}
\put(16.5,44){$\beta_2$}
\put(27.5,44){$\beta_3$}
\put(60.5,44){$\beta_4$}
\put(84.2,44){$\beta_5$}
\put(88,44){$\beta_6$}
\put(44,1){mode coefficient, $j$}
\put(3,9.5){$|{\bf \hat a}_j(t_3)|$}
\put(3,23){$|{\bf \hat a}_j(t_2)|$}
\put(3,37){$|{\bf \hat a}_j(t_1)|$}
\end{overpic}
\end{center}
\vspace{-.1in}
\caption{\label{fig:cqgle3} Full evolution dynamics (a) and the low-rank POD dynamic reconstruction using
compressive sensing and Galerkin projection (b).  The black lines on
the left panel at $t=25, 125$ and 225 represent the
sampling times while the three black circles represent the
three sparse measurement locations.  From the three samples,
the right panel is reconstructed by identifying the correct
POD modes and using Galerkin projection.  The bottom panel
shows the modal coefficient vector ${\bf \hat a}$ evaluated via convex $\ell_1$ optimization for the
three sampling times.  Correct identification is achieved  of the $\beta_1$ regime
at $t_1=25$ (circles, $\circ$), of the $\beta_3$ regime at $t_2=125$ (triangles, $\triangle$), and of
the $\beta_5$ regime at $t_3=225$ (plusses, $+$).  The $|{\bf \hat a}_j|$ are color coded
according to the library elements depicted in Figure~\ref{fig:cqgle2}.  For ease of
viewing, the different $\beta_j$ regimes are separated by shaded/non-shaded
regions and are further identified at the top of panel (c).}
\end{figure}

To highlight the role of compressive sensing in identification and reconstruction for dynamical systems,  we allow the bifurcation parameter $\beta=\beta(t)$ to vary in time so that the dynamics switch between attractors as $\beta$ changes.  
Consider an example where $\beta=\beta_1$ for $t\in[0,100)$, $\beta=\beta_3$ for $t\in[100,200)$ and $\beta=\beta_5$ for $t\in[200,300]$.  
The evolution dynamics for this case are illustrated in Figure~\ref{fig:cqgle3} (a).

We measure the state at either three ($x_1$-$x_3$) or five ($x_1$-$x_5$) locations $x_1=0, x_2=0.7, x_3=1.4, x_4=1.8,x_5=2.2$ (shown at the bottom of Figure~\ref{fig:cqgle2}) taking data only at the times $t_1=25, t_2=125$ and $t_3=225$; these times are chosen $25$ units after the bifurcation value switches so that transients have decayed.  
At each instance, we take sparse measurements and perform classification, projection, and forward simulation (Galerkin reconstruction), while working exclusively in the low-dimensional POD library.  
The procedure is as follows:
\begin{itemize}
\item[(i)] \textbf{[classification]} from a sparse set of measurements (three or five), the modes corresponding to the specific $\beta_j$ are identified and extracted, 
\item[(ii)] \textbf{[projection]} the sparse measurements are projected, through a standard pseudo-inverse operation, onto the modes ${\bf\Psi}_j$ for the particular parameter $\beta_j$ to determine initial values of $a_n$, 
\item[(iii)] \textbf{[reconstruction]} the extracted library modes are evolved according the the POD-Galerkin projection technique by using the spatial modes from the library ${\bf \Psi}_j$ in conjunction with their time dynamics $a_n(t)$~\cite{HLBR_turb}.
\end{itemize}
Figure~\ref{fig:cqgle3} (b) shows the resulting dynamic reconstruction, and Figure~\ref{fig:cqgle3} (c) shows the coefficients of the specific sparse vector $\bf{\hat{a}}$ identified at each time: $t_1$, $t_2$, and $t_3$.  
Indeed, the proposed algorithm using only three measurements reproduces the dynamics with remarkable success.  
Note that the number of measurements $m=3$, number of library elements $p=24$, and the original size of the system $n=1024$ are ordered so that $m\ll p \ll n$.  
Consequently, the matrix ${\bf\Phi \Psi}$ in the under-determined system \eqref{eq:sparse3} is a $3\times 24$ matrix, yielding an efficient $\ell_1$ convex optimization problem for the sparse identification.

For this example, we choose effective sampling locations based on the library modes of ${\bf \Psi}$.  
If poor choices are made, i.e. not aligning the sensors with maxima and minima observed in the POD library modes~\cite{Yildirim:2009}, then the dominant modes are often misidentified.
This sensitivity to sensor location suggests that sensor placement should be carefully considered.  
Moreover, it is assumed that measurements of the system are perfect.  
However, noise is inherent in the detector and/or model.

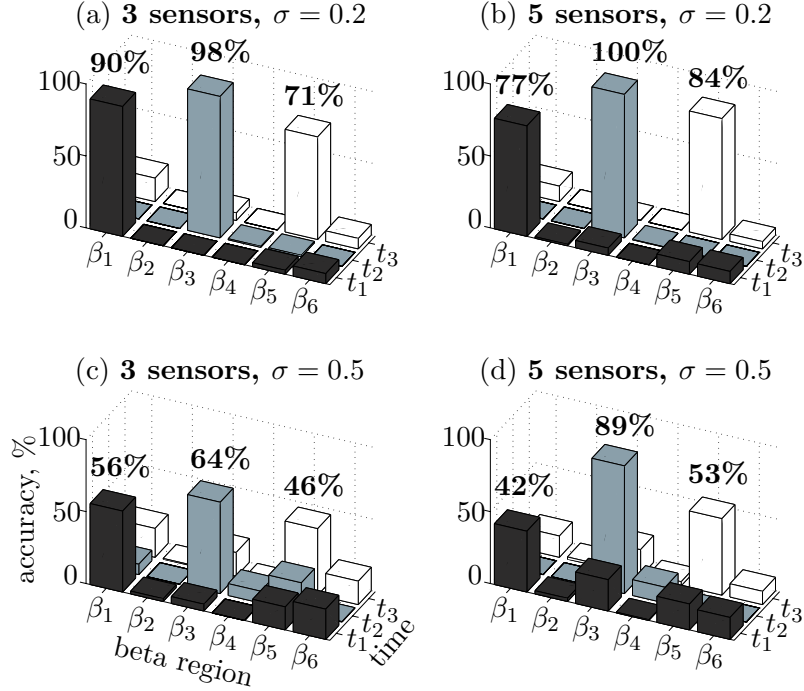
\begin{figure}
\begin{center}
\begin{tabular}{cc}
\begin{overpic}[width=.3\textwidth]{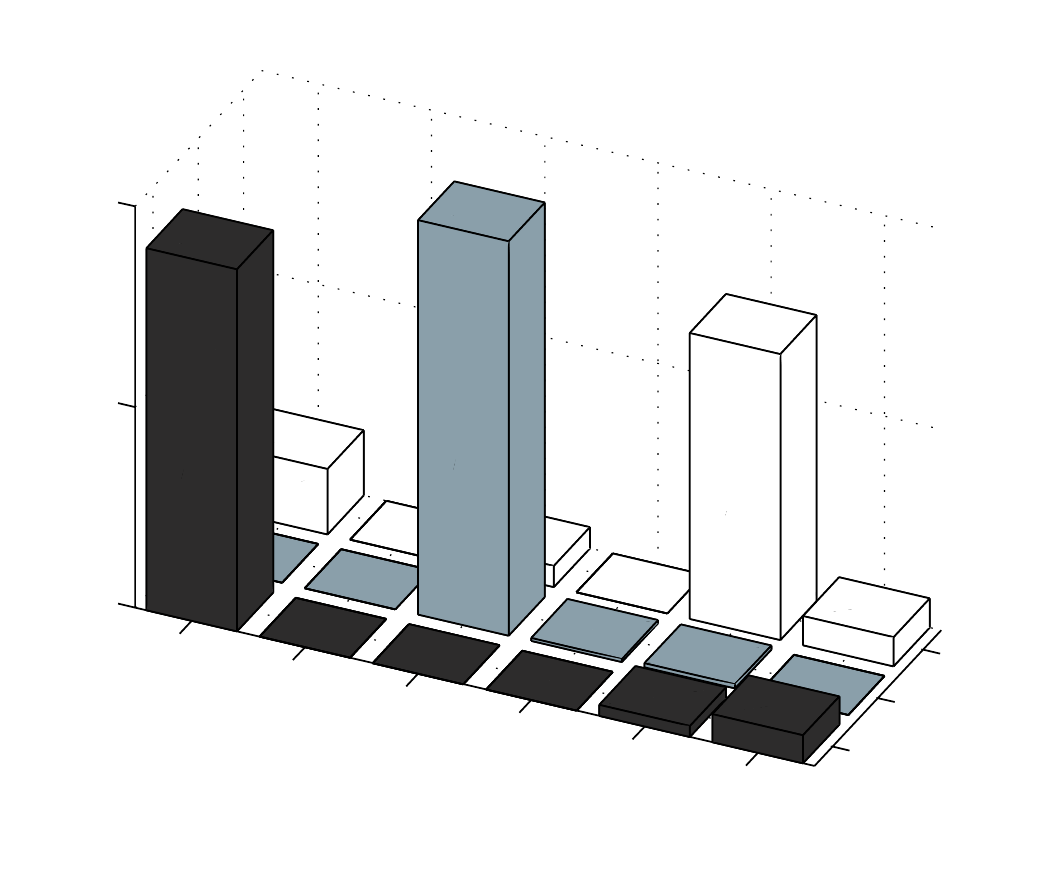}
\put(10,80){(a) \bf{3 sensors, $\sigma=0.2$}}
\put(14,67){\bf 90\%}
\put(41,70){\bf 98\%}
\put(66,58){\bf 71\%}
\put(7,25){0}
\put(4,42){50}
\put(0,62){100}
\put(13,16){$\beta_1$}
\put(24,13.5){$\beta_2$}
\put(35,11){$\beta_3$}
\put(46,8.5){$\beta_4$}
\put(57,6){$\beta_5$}
\put(68,3.5){$\beta_6$}
\put(83,8){$t_1$}
\put(87,12){$t_2$}
\put(91,16){$t_3$}
\end{overpic} & 
\begin{overpic}[width=.3\textwidth]{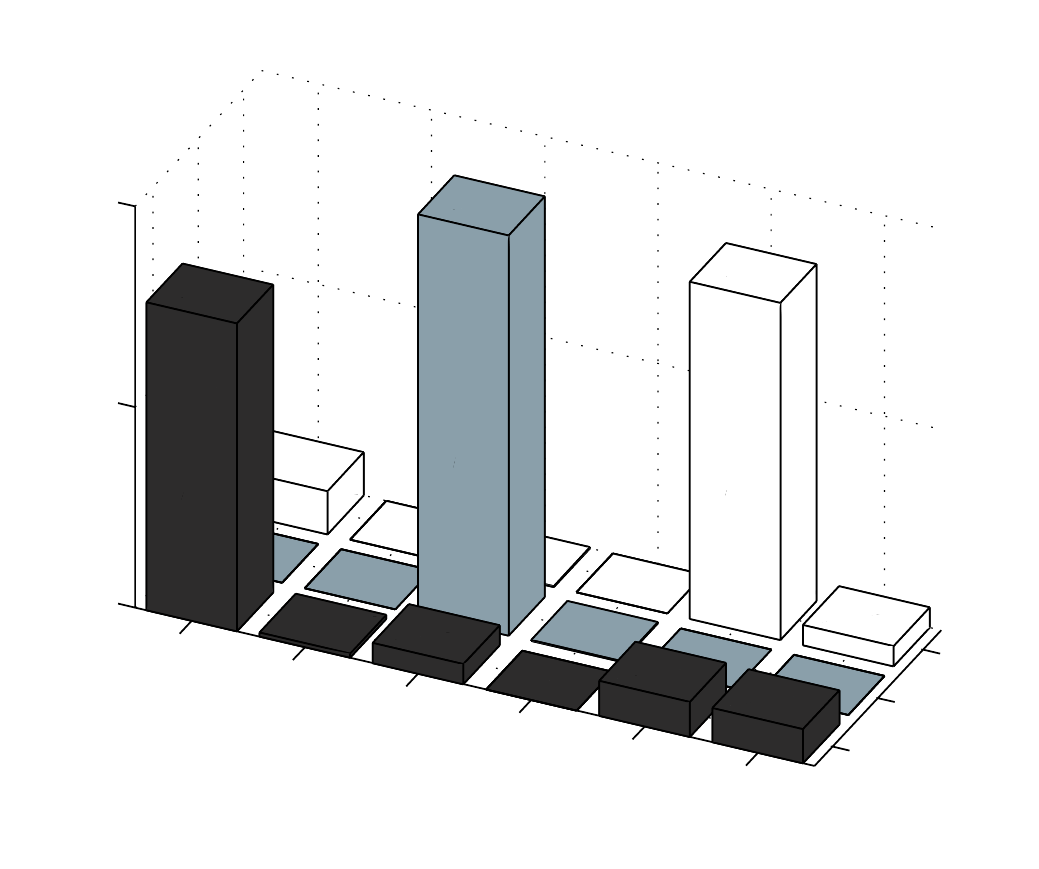} 
\put(10,80){(b) \bf{5 sensors, $\sigma=0.2$}}
\put(14,61){\bf 77\%}
\put(40,70){\bf 100\%}
\put(66,63){\bf 84\%}
\put(7,25){0}
\put(4,42){50}
\put(0,62){100}
\put(13,16){$\beta_1$}
\put(24,13.5){$\beta_2$}
\put(35,11){$\beta_3$}
\put(46,8.5){$\beta_4$}
\put(57,6){$\beta_5$}
\put(68,3.5){$\beta_6$}
\put(83,8){$t_1$}
\put(87,12){$t_2$}
\put(91,16){$t_3$}
\end{overpic}\\
&\\
\begin{overpic}[width=.3\textwidth]{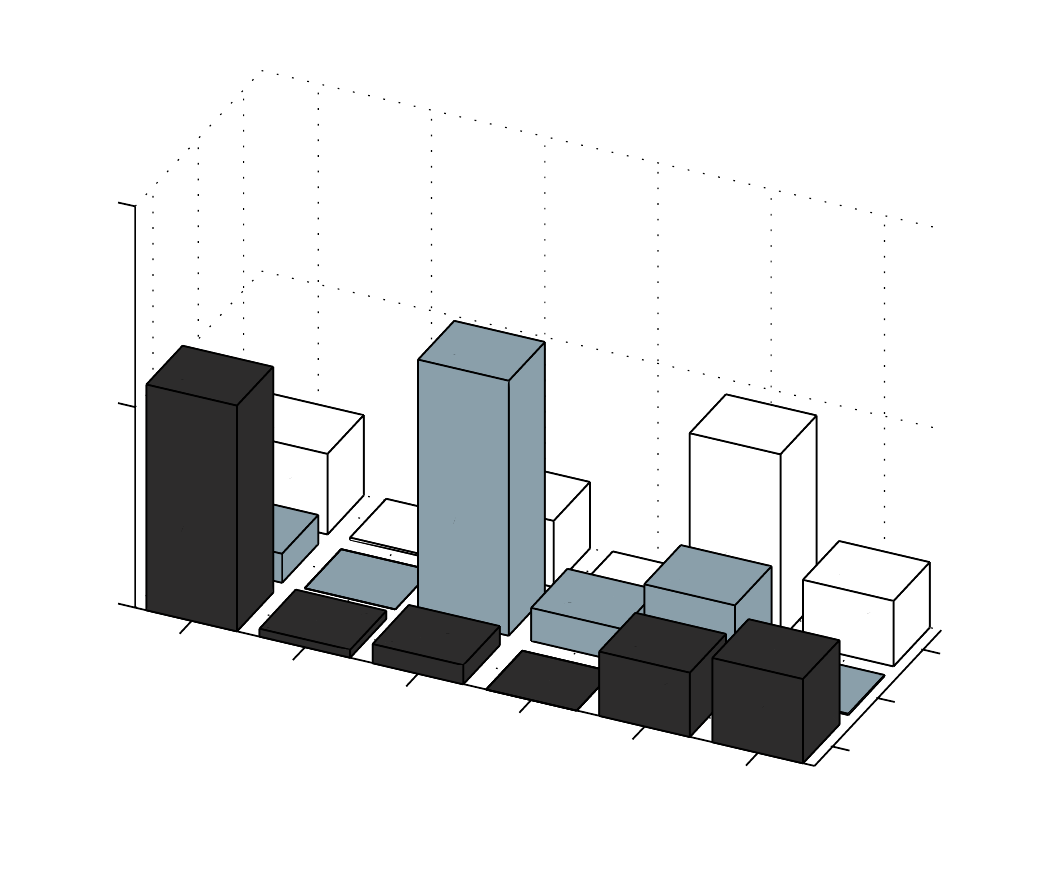} 
\put(10,80){(c) \bf{3 sensors, $\sigma=0.5$}}
\put(20,6){\begin{rotate}{-13}beta region\end{rotate}}
\put(-1,25){\begin{rotate}{90}accuracy, \%\end{rotate}}
\put(92,2){\begin{rotate}{50}time\end{rotate}}
\put(14,54){\bf 56\%}
\put(41,56){\bf 64\%}
\put(66,49){\bf 46\%}
\put(7,25){0}
\put(4,42){50}
\put(0,62){100}
\put(13,16){$\beta_1$}
\put(24,13.5){$\beta_2$}
\put(35,11){$\beta_3$}
\put(46,8.5){$\beta_4$}
\put(57,6){$\beta_5$}
\put(68,3.5){$\beta_6$}
\put(83,8){$t_1$}
\put(87,12){$t_2$}
\put(91,16){$t_3$}
\end{overpic}& 
\begin{overpic}[width=.3\textwidth]{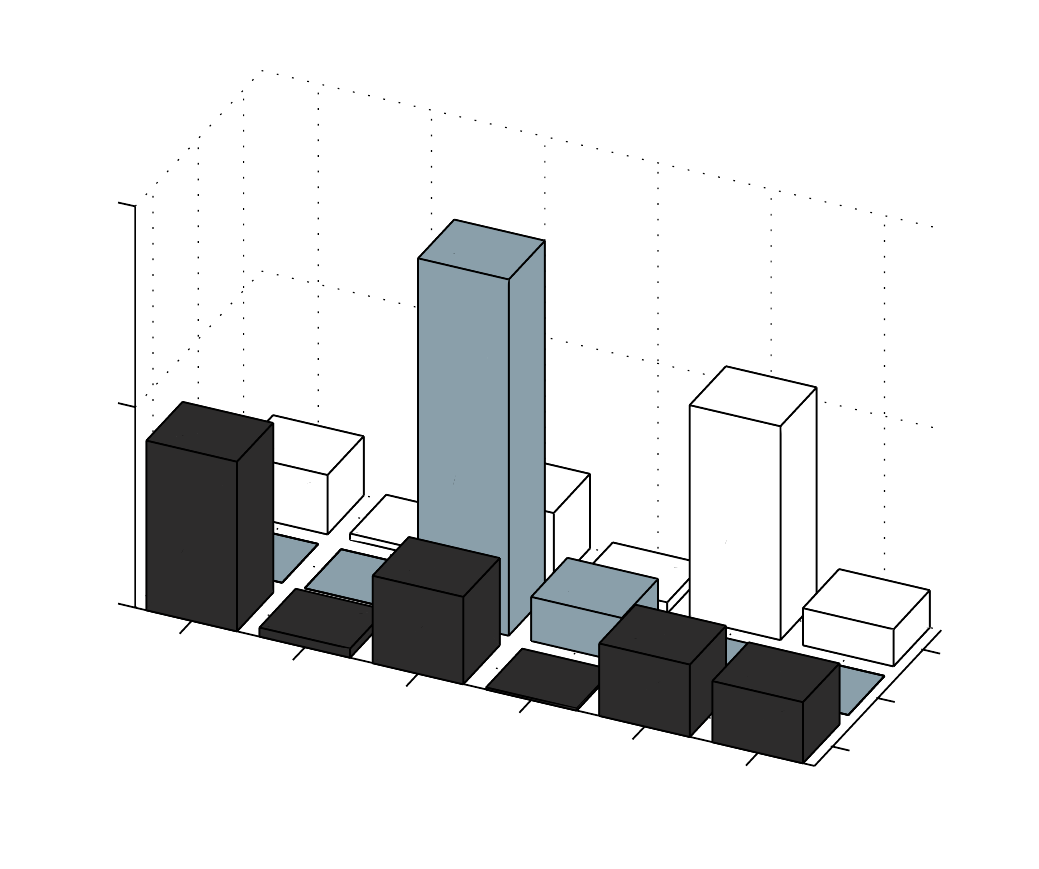} 
\put(10,80){(d) \bf{5 sensors, $\sigma=0.5$}}
\put(14,49){\bf 42\%}
\put(41,65){\bf 89\%}
\put(66,52){\bf 53\%}
\put(7,25){0}
\put(4,42){50}
\put(0,62){100}
\put(13,16){$\beta_1$}
\put(24,13.5){$\beta_2$}
\put(35,11){$\beta_3$}
\put(46,8.5){$\beta_4$}
\put(57,6){$\beta_5$}
\put(68,3.5){$\beta_6$}
\put(83,8){$t_1$}
\put(87,12){$t_2$}
\put(91,16){$t_3$}
\end{overpic} 
\end{tabular}
\caption{\label{fig:cqgle_ml_6} 
Accuracy of region classification for the system described in Figure~\ref{fig:cqgle3}.  
The bifurcation parameter $\beta$ switches from $\beta_1$ at $t_1=25$ to $\beta_3$ at $t_2=125$ to $\beta_5$ at $t_3=225$.  
The bar charts illustrate which bifurcation regime $\beta_1$-$\beta_6$ is classified from sparse measurements by the $\ell_1$-minimization procedure described above.  
Three or five sensors are considered under moderate $\sigma=0.2$ and strong $\sigma=0.5$ error measurements using 400 realizations.   
More sensors improve the region identification performance for regions $\beta_3$ and $\beta_5$, but decrease performance for $\beta_1$.}
\end{center}
\end{figure}

To quantify the impact of noise on the classification and reconstruction, \eqref{eq:sparse2} is modified to {${\bf \hat{U}}= {\bf\Phi} {\bf U} + {\mathcal{N}}(0,\sigma^2)$} where Gaussian distributed, white-noise error $\mathcal{N}$ with variance $\sigma^2$ is added to account for measurement error.  
Figure~\ref{fig:cqgle_ml_6} shows statistical results of 400 trials using three or five sensors for noise strength $\sigma=0.2$ or 0.5.  
With moderate noise ($\sigma=0.2$), both the three and five sensor scenarios identify the correct regime quite accurately.  
For stronger noise ($\sigma=0.5$), both three and five sensors lose a great deal of accuracy in the identification process.
It is also observed that more sensors actually hinders the evaluation of the $\beta_1$ parameter regime, although the $\beta_3$ and $\beta_5$ cases improve.  
Indeed, numerical simulations indicate that the three sensors placed at $x=0, 0.7$, and 1.4 are robust and are not easily improved on by varying placement or quantity.  
Further study is needed to determine optimal sensor location.

These results suggest that multiple samplings in time can be used to reach a statistical conclusion about the correct parameter regime, thus avoiding mis-identification.  
For example, we already wait 25 time units after $\beta$ switches to take sparse measurements, so that transients decay.  
If, instead of sampling a single time unit at $t=25$, we accumulate information over 5-10 time units, the effect of sensor noise is significantly attenuated.  

We also investigate the least-squares estimate $\bf\tilde a$ for the mode amplitudes based on the 3-sensor and 5-sensor configurations.  
In every single case, for no noise, as well as for noise levels $\sigma=0.2$ and $\sigma=0.5$, the least-squares solution $\bf \tilde a$ results in the misidentification of $\beta_1$ and $\beta_5$ regimes, instead identifying the incorrect $\beta_3$ regime.  
The collapse of $\ell_2$ minimization for identifying the bifurcation parameter regime highlights the success of the sparse sampling strategy, centered around $\ell_1$-minimization.

\section{Discussion}\label{sec:discussion}
In conclusion, we advocate a general theoretical framework for complex systems whereby low-rank structures are represented by the $\ell_2$-optimal proper orthogonal decomposition, and then identified from limited noisy measurements using the sparsity promoting $\ell_1$ norm and the compressive sensing architecture.  
The strategy for building a modal library by concatenating truncated POD libraries across a range of relevant bifurcation parameters may be viewed as a simple machine learning implementation.  
The resulting modal library is a natural sparse basis for the application of compressive sensing.  After the expensive one-time library-building procedure, accurate identification, projection, and reconstruction may be performed entirely in a low-dimensional framework.  

To our knowledge, these results are the first to combine even simple machine learning concepts and compressive sensing to complex systems for both
\begin{itemize}
\item[(i)] correctly identifying the dynamical parameter regime, and
\item[(ii)] reconstructing the associated low-rank dynamics.
\end{itemize}
Pairing a low-dimensional learned library, in which the dynamics have a sparse representation, with compressive sensing provides a powerful new architecture for studying dynamical systems that exhibit coherent behavior.

With three sensors, it is possible to accurately classify bifurcation regime, reconstruct the low-dimensional content, and simulate the Galerkin projected dynamics of the complex Ginzburg Landau equation.  
In addition, we investigate the performance of compressed sensing with the addition of sensor noise and the addition of more sensors.  
For moderate noise levels, the method accurately classifies the correct dynamic regime, although performance drops for larger noise values.  
The addition of more sensors does not significantly improve performance, although the sensor placement was not exhaustive.  
In contrast, classification based on least-squares fails to identify $\beta_1$ and $\beta_5$ regions for all noise levels, on every trial.

There are a number of important directions that arise from this work.  
The library building procedure discussed in Figure~\ref{fig:L2schematic} is quite general, and it will be interesting to investigate additional library building techniques and machine learning strategies.  
For example, is it possible to remove features that are common to all of the dynamic regimes to enhance contrast between categories in the $\ell_1$ classification step?  
It will also be interesting to investigate optimal sensor placement based on the principle of maximizing incoherence with respect to the overcomplete basis.  
Finally, it may be possible to use coherence between each pair of local bases $({\bf \Psi}_i,{\bf\Psi}_j)$ as a means to construct an induced metric on the space of bifurcation parameters.  
This may facilitate the accurate categorization of dynamic regimes that have not been directly explored in the training step.  
Finally, the procedure above is promising for use with data assimilation techniques, which typically incorporate new measurements using least-squares fitting ($\ell_2$).

As these directions unfold, we believe that the combination of $\ell_2$ low-rank representations and $\ell_1$ sparse sampling will enable efficient characterization and manipulation of low-rank dynamical systems.  
The ultimate goal is to always work in a measurement space with dimension on the order of the underlying low-dimensional attractor.

\section*{Acknowledgment}
We are grateful for discussions with Bingni W. Brunton, Joshua L. Proctor, and Xing Fu.
J. N. Kutz acknowledges support
from the National Science Foundation (DMS-1007621) and the
U.S. Air Force Office of Scientific Research (FA9550-09-0174).  

\bibliographystyle{plain}
\bibliography{jabbrv.bib,references.bib}
\end{document}